\documentclass[a4paper,reqno,10pt]{amsart}

\usepackage{graphicx}
\usepackage{amsmath,mathtools,latexsym,amssymb,amsthm}
\allowdisplaybreaks
\usepackage[nospace,noadjust]{cite}
\usepackage{stmaryrd}
\usepackage{cancel}
\usepackage[all,cmtip]{xy}
\usepackage{enumerate}
\usepackage{setspace}
\usepackage[mathscr]{euscript}

\theoremstyle{plain}
\newtheorem{thm}{Theorem}[section]
\newtheorem{prop}[thm]{Proposition}
\newtheorem{lemma}[thm]{Lemma}

\theoremstyle{definition}
\newtheorem{defn}[thm]{Definition}
\newtheorem{remark}[thm]{Remark}
\newtheorem{exam}[thm]{Example}

\newcommand{\Lt}[1]{\mathbin{\triangleleft_{#1}}}
\newcommand{\Rt}[1]{\mathbin{\triangleright_{#1}}}
\newcommand{\g}{\mathfrak{a}}
\newcommand{\h}{\mathfrak{b}}
\newcommand{\pL}{\mathrm{pL}}
\newcommand{\Leib}{\mathrm{Leib}}
\newcommand{\Tot}{\mathrm{Tot}}
\newcommand{\kk}{{\bf k}}

\usepackage{hyperref}
\hypersetup{colorlinks, citecolor=blue, filecolor=black, linkcolor=blue, urlcolor=black}

\begin{document}

\title[Weighted \texorpdfstring{$\mathcal{O}$}{O}-operators on pre-Leibniz algebras]{Weighted \texorpdfstring{$\mathcal{O}$}{O}-operators on pre-Leibniz algebras and post-pre-Leibniz structures}

\author{Sami Benabdelhafidh}
\address{University of Sfax, Faculty of Sciences of Sfax, BP 1171, 3038 Sfax, Tunisia}
\renewcommand{\emailaddrname}{Email addresses}
\email{abdelhafidhsami41@gmail.com}
%samibenabdelhafidh95@gmail.com}

\subjclass[2020]{17A32, 17A36, 17B56, 16S80.}
\keywords{Pre-Leibniz algebra, Weighted $\mathcal{O}$-operator, Rota-Baxter operator, Cohomology, Deformation, Post-pre-Leibniz algebra.}

\begin{abstract}
%Pre-Leibniz algebras (also called Leibniz-dendriform algebras) are a non-skewsymmetric analogue of pre-Lie algebras: the sum of their two operations is a Leibniz bracket. 
In this paper, we introduce weighted $\mathcal{O}$-operators (also called weighted relative Rota-Baxter operators) on pre-Leibniz algebras. We construct a differential graded Lie algebra whose Maurer-Cartan elements are exactly these operators. This allows us to define their cohomology as that of the induced pre-Leibniz algebra with coefficients in a suitable representation. As an application, we study formal and finite order deformations of these operators. We show that this theory splits the corresponding theory for weighted $\mathcal{O}$-operators on Leibniz algebras. Finally, we introduce post-pre-Leibniz algebras as the algebraic structures underlying weighted $\mathcal{O}$-operators on pre-Leibniz algebras.
\end{abstract}

\maketitle

\tableofcontents

\thispagestyle{empty}

\section{Introduction}\label{sec-1}

Rota-Baxter operators first appeared in the work of Baxter on fluctuation theory in probability \cite{baxter}, and were later studied by Rota \cite{rota} and Cartier \cite{cartier}. They have applications in the splitting of operads \cite{bai-operad}, in the Connes-Kreimer approach to renormalization \cite{connes} and in noncommutative Poisson geometry \cite{uchino}. An $\mathcal{O}$-operator (or relative Rota-Baxter operator), introduced by Kupershmidt \cite{kuper}, is a generalization of a Rota-Baxter operator in the presence of a representation. $\mathcal{O}$-operators of weight $\lambda$ were introduced by Bai, Guo and Ni \cite{bgn, bgn2} in their study of the classical Yang-Baxter equation and post-Lie algebras. It is well known that an operator of weight zero splits an algebra into a pre-algebra, whereas an operator of nonzero weight gives a post-algebra \cite{bai-operad, ebrahimi}. Post-Lie algebras were introduced by Vallette \cite{vallette}.

\medskip

The deformation and cohomology theories were first developed for operators of weight zero \cite{tang, das-rota, lst}. More recently, operators of arbitrary weight have received a lot of attention. Das studied the cohomology and deformations of weighted Rota-Baxter operators \cite{das-weighted}, and of weighted relative Rota-Baxter operators on Leibniz algebras together with the associated post-Leibniz algebras \cite{das-wrb}. In \cite{bbt-jj}, the cohomology of modified $\lambda$-differential Jacobi-Jordan algebras was developed, together with its applications. The cohomology of Rota-Baxter $3$-Lie algebras of arbitrary weight was studied in \cite{gqwz}. Hou, Sheng and Zhou \cite{hsz} studied relative Rota-Baxter operators of nonzero weight on $3$-Lie algebras, their deformations and the induced $3$-post-Lie algebras. Modified Rota-Baxter Lie-Yamaguti algebras and their cohomology and deformations were considered in \cite{bb-lya, tg-lya}. The cohomology and deformations of morphisms of modified $\lambda$-differential associative algebras were studied in \cite{b-mor}. Weighted $\mathcal{O}$-operators on Malcev algebras were introduced in \cite{hmns}, together with post-Malcev algebras. Weighted $\mathcal{O}$-operators on Hom-Lie triple systems were considered in \cite{teng-hlts}. On the group side, Rota-Baxter operators of weight $1$ on Lie groups were introduced in \cite{gls}; they are related to post-groups \cite{bgst}, whose differentiation gives post-Lie algebras.

\medskip

Leibniz algebras are a non-skewsymmetric version of Lie algebras \cite{bloh, loday-une}. Their universal enveloping algebras and their cohomology were introduced by Loday and Pirashvili \cite{loday-pira}. Relative Rota-Baxter operators of weight $0$ on Leibniz algebras were introduced by Tang and Sheng \cite{sheng-pl}, and their deformations and related structures were studied in \cite{sheng-defr, sun-jing}. These operators give rise to pre-Leibniz algebras, introduced in \cite{sheng-pl} under the name of Leibniz-dendriform algebras, which split Leibniz algebras just as pre-Lie algebras split Lie algebras. Pre-Leibniz algebras were then studied in detail by Das \cite{das-pre}: he described them as Maurer-Cartan elements of a graded Lie algebra \cite{bala-operad} and developed their cohomology, which splits the Loday-Pirashvili cohomology.

\medskip

In this paper, we extend the results of \cite{das-wrb} to pre-Leibniz algebras, using the methods of \cite{das-pre}. Let $\g$ and $\h$ be pre-Leibniz algebras such that $\h$ is a pre-Leibniz $\g$-representation. A linear map $T : \h \rightarrow \g$ is called a $\lambda$-weighted $\mathcal{O}$-operator if
\[
T(u) \ast T(v) = T \big( T(u) \ast v + u \ast T(v) + \lambda \, u \ast v \big), \qquad \forall u, v \in \h, ~ \ast \in \{ \triangleleft, \triangleright \}.
\]
We characterize such operators by their graphs in the weighted semidirect product $\g \ltimes_\lambda \h$, and we show that $T$ induces a new pre-Leibniz algebra $\h_T$. Using the graded Lie algebra of \cite{das-pre} and Voronov's derived brackets \cite{voro}, we construct a differential graded Lie algebra whose Maurer-Cartan elements are the $\lambda$-weighted $\mathcal{O}$-operators (Theorem \ref{mc-thm}). We then make $\g$ into a representation of $\h_T$ and define the cohomology of $T$ as the pre-Leibniz cohomology of $\h_T$ with coefficients in $\g$. As an application, we study formal and finite order deformations of $T$. We give a sufficient condition for rigidity, and we show that a finite order deformation extends if and only if its obstruction class vanishes.

\medskip

We also introduce, following \cite{vallette, hsz, hmns}, the notion of post-pre-Leibniz algebras. They have six operations satisfying twenty-one identities. Every $\lambda$-weighted $\mathcal{O}$-operator induces a post-pre-Leibniz algebra and, conversely, every post-pre-Leibniz algebra comes from a $\lambda$-weighted $\mathcal{O}$-operator. Adding the operations of a post-pre-Leibniz algebra, we obtain a post-Leibniz algebra in the sense of \cite{das-wrb}. Thus we have the following commutative diagram:
\[
\xymatrix@C=4pc@R=1.5pc{
\text{post-pre-Leibniz algebras} \ar[r] \ar[d] & \text{pre-Leibniz algebras} \ar[d]\\
\text{post-Leibniz algebras} \ar[r] & \text{Leibniz algebras}
}
\]

The paper is organized as follows. In Section \ref{sec-2}, we recall
some background
about pre-Leibniz algebras. In Section \ref{sec-3}, we introduce $\lambda$-weighted $\mathcal{O}$-operators on pre-Leibniz algebras, give some examples and describe the induced pre-Leibniz algebra. In Section \ref{sec-4}, we construct the differential graded Lie algebra whose Maurer-Cartan elements are $\lambda$-weighted $\mathcal{O}$-operators and we define their cohomology. In Section \ref{sec-5}, we study formal and finite order deformations, rigidity and obstructions. Finally, in Section \ref{sec-6}, we introduce post-pre-Leibniz algebras and study their relation with $\lambda$-weighted $\mathcal{O}$-operators and post-Leibniz algebras.

\medskip

Throughout the paper, all vector spaces, linear and multilinear maps and tensor products are over a field $\kk$ of characteristic $0$, and $\lambda \in \kk$ is a fixed scalar.

\section{Background on pre-Leibniz algebras}\label{sec-2}

In this section, we recall Leibniz algebras, pre-Leibniz algebras and their cohomology. Our main references are \cite{loday-pira, bala-operad, das-pre}.

\begin{defn}
A Leibniz algebra is a vector space $\mathfrak{g}$ with a bilinear bracket $[\cdot, \cdot]$ satisfying
\[
[x, [y,z]] = [[x,y],z] + [y,[x,z]], \qquad \forall x, y, z \in \mathfrak{g}.
\]
A representation of $\mathfrak{g}$ is a vector space $V$ with bilinear maps $\rho^L : \mathfrak{g} \otimes V \rightarrow V$ and $\rho^R : V \otimes \mathfrak{g} \rightarrow V$ such that $\mathfrak{g} \oplus V$ is a Leibniz algebra with the bracket
\[
[(x,u), (y,v)] = \big( [x,y], ~ \rho^L (x,v) + \rho^R (u,y) \big), \qquad \forall x, y \in \mathfrak{g}, ~ u, v \in V.
\]
\end{defn}

The Loday-Pirashvili cochain complex of $\mathfrak{g}$ with coefficients in $V$ is given by $C^n_\Leib (\mathfrak{g}, V) = \mathrm{Hom}(\mathfrak{g}^{\otimes n}, V)$ and
\begin{equation}\label{lp-diff}
\begin{split}
(\delta_{\Leib} f ) (x_1, \ldots, x_{n+1}) 
&= \sum_{i=1}^n (-1)^{i+1} \rho^L \big( x_i, f(x_1, \ldots, \widehat{x_i}, \ldots, x_{n+1}) \big) \\
&{}+ (-1)^{n+1} \rho^R \big( f (x_1, \ldots, x_n), x_{n+1} \big) \\
&{}+ \sum_{1 \leq i < j \leq n+1} (-1)^i  f \big( x_1, \ldots, \widehat{x_i}, \ldots, x_{j-1}, [x_i, x_j], x_{j+1}, \ldots, x_{n+1} \big).
\end{split}
\end{equation}
Its cohomology is denoted by $H^\ast_\Leib (\mathfrak{g}, V)$ \cite{loday-pira}.

\begin{defn}\label{defn-pre-leibniz}\cite{sheng-pl, das-pre}
A pre-Leibniz algebra is a triple $(\g, \triangleleft, \triangleright)$ consisting of a vector space $\g$ with two bilinear operations $\triangleleft, \triangleright : \g \otimes \g \rightarrow \g$ satisfying, $\forall x, y, z \in \g$,
\begin{align}
x \triangleleft ( y \triangleleft z + y \triangleright z ) &= ( x \triangleleft y ) \triangleleft z + y \triangleright (x \triangleleft z), \label{pre-l1}\\
x \triangleright ( y \triangleleft z ) &= ( x \triangleright y ) \triangleleft z + y \triangleleft ( x \triangleleft z + x \triangleright z),\label{pre-l2}\\
x \triangleright ( y \triangleright z) &= ( x  \triangleleft y + x \triangleright y) \triangleright z + y \triangleright ( x \triangleright z). \label{pre-l3}
\end{align}
A morphism of pre-Leibniz algebras and a derivation of $\g$ are linear maps $\phi$ and $D$ satisfying respectively
\[
\phi (x \ast y) = \phi(x) \ast' \phi(y) \quad \text{and} \quad D (x \ast y) = D(x) \ast y + x \ast D(y), \qquad \forall x, y \in \g, ~ \ast \in \{ \triangleleft, \triangleright \}.
\]
\end{defn}

Throughout the paper, we write $\ast \in \{ \triangleleft, \triangleright \}$ to state two identities at once. If $(\g, \triangleleft, \triangleright)$ is a pre-Leibniz algebra, then $\g_\Tot = (\g, [\cdot, \cdot]_{\triangleleft, \triangleright})$ is a Leibniz algebra, called the total Leibniz algebra \cite{sheng-pl, das-pre}, where
\[
[x, y]_{\triangleleft, \triangleright} = x \triangleleft y + x \triangleright y, \qquad \forall x, y \in \g.
\] A pre-Leibniz algebra with $x \triangleleft y = - y \triangleright x$ is the same as a pre-Lie algebra $(\g, \triangleleft)$ \cite{das-pre}.

\begin{exam}\label{exam-pre}
(i) If $(\mathfrak{g}, [\cdot, \cdot])$ is a Leibniz algebra, then $(\mathfrak{g}, \triangleleft = 0, \triangleright = [\cdot, \cdot])$ is a pre-Leibniz algebra: (\ref{pre-l1}) and (\ref{pre-l2}) are trivial and (\ref{pre-l3}) is the Leibniz identity.

(ii) If $(A, \cdot)$ is an associative algebra, then $A$ is a pre-Leibniz algebra with
\[
x \triangleleft y = x \cdot y, \qquad x \triangleright y = - y \cdot x, \qquad \forall x, y \in A.
\]
For instance, (\ref{pre-l1}) reads $x(yz - zy) = (xy)z - (xz)y$.
\end{exam}

\begin{defn}\label{defn-rep}
A representation of a pre-Leibniz algebra $(\g, \triangleleft, \triangleright)$ is a vector space $M$ with four bilinear maps $\triangleleft^L, \triangleright^L : \g \otimes M \rightarrow M$ and $\triangleleft^R, \triangleright^R : M \otimes \g \rightarrow M$ such that the direct sum $\g \oplus M$ is a pre-Leibniz algebra (the semidirect product $\g \ltimes M$) with the operations
\[
(x, u) \ast (y, v) = \big( x \ast y, ~ x \ast^L v + u \ast^R y \big), \qquad \forall x, y \in \g, ~ u, v \in M, ~ \ast \in \{ \triangleleft, \triangleright \}.
\]
\end{defn}

Equivalently, the nine identities obtained from (\ref{pre-l1})-(\ref{pre-l3}) by replacing exactly one of $x, y, z$ by an element of $M$ hold. To lighten the notation, we often write $x \ast u$ and $u \ast x$ instead of $x \ast^L u$ and $u \ast^R x$. Any pre-Leibniz algebra is a representation of itself (the adjoint representation). If $M$ is a representation of $\g$, then $M_\Tot = (M, \rho^L, \rho^R)$ is a representation of the Leibniz algebra $\g_\Tot$ \cite{das-pre}, where
\[
\rho^L = \triangleleft^L + \triangleright^L, \qquad \rho^R = \triangleleft^R + \triangleright^R.
\]

\medskip

\noindent The graded Lie algebra. Pre-Leibniz algebras can be described as Maurer-Cartan elements of a graded Lie algebra, which we now recall from \cite{das-pre}. Let $C_n = \{ [1], \ldots, [n] \}$ and let $\kk[C_n]$ be the vector space spanned by $C_n$. Let $m, n \geq 1$, $1 \leq i \leq m$ and $\sigma \in Sh(i-1, n-1)$. We put the elements of $C_{m+n-1}$ into $m$ ordered blocks as follows:
\[
\{ [\sigma(1)] \}, \ldots, \{ [\sigma(i-1)] \}, ~ \{ [\sigma(i)], \ldots, [\sigma(i+n-2)], [i+n-1] \}, ~ \{ [i+n] \}, \ldots, \{ [m+n-1] \}.
\]
Define maps $R^\sigma_{m;i,n} : C_{m+n-1} \rightarrow C_m$ and $S^\sigma_{m;i,n} : C_{m+n-1} \rightarrow \kk[C_n]$ by $R^\sigma_{m;i,n} ([r]) = [k]$ if $[r]$ lies in the $k$-th block, and
\[
S^\sigma_{m;i,n} ([r]) = \begin{cases} [1] + \cdots + [n] & \text{if } [r] \text{ is not in the } i\text{-th block},\\
[j] & \text{if } [r] \text{ is the } j\text{-th element of the } i\text{-th block}. \end{cases}
\]
For a vector space $V$, let $C^n_\pL (V, V) = \mathrm{Hom}(\kk[C_n] \otimes V^{\otimes n}, V)$. For $f \in C^m_\pL (V,V)$, $g \in C^n_\pL (V,V)$ and $1 \leq i \leq m$, define $f \diamond_i g \in C^{m+n-1}_\pL (V,V)$ by
\begin{multline*}
(f \diamond_i g)([r]; x_1, \ldots, x_{m+n-1}) = \sum_{\sigma \in Sh(i-1, n-1)} (-1)^\sigma f \Big( R^\sigma_{m;i,n}[r]; x_{\sigma(1)}, \ldots, x_{\sigma(i-1)}, \\
g \big( S^\sigma_{m;i,n}[r]; x_{\sigma(i)}, \ldots, x_{\sigma(i+n-2)}, x_{i+n-1} \big), x_{i+n}, \ldots, x_{m+n-1} \Big),
\end{multline*}
and
\begin{align}\label{pl-bracket}
\llbracket f, g \rrbracket_\pL = \sum_{i=1}^m (-1)^{(i-1)(n-1)} f \diamond_i g - (-1)^{(m-1)(n-1)} \sum_{i=1}^n (-1)^{(i-1)(m-1)} g \diamond_i f.
\end{align}

In particular, let $m = n = 2$. For $i = 1$, the only shuffle is the identity and the blocks are $\{ [1], [2] \}, \{ [3] \}$; for $i = 2$, the shuffles are $\mathrm{id}$ and the transposition $(12)$ of $Sh(1,1)$, with blocks $\{ [1] \}, \{ [2], [3] \}$ and $\{ [2] \}, \{ [1], [3] \}$ respectively. Therefore
\begin{align*}
(f \diamond_1 g)([1]; x_1, x_2, x_3) &= f \big( [1]; g([1]; x_1, x_2), x_3 \big), \\
(f \diamond_1 g)([2]; x_1, x_2, x_3) &= f \big( [1]; g([2]; x_1, x_2), x_3 \big), \\
(f \diamond_1 g)([3]; x_1, x_2, x_3) &= f \big( [2]; g([1] + [2]; x_1, x_2), x_3 \big), \\
(f \diamond_2 g)([1]; x_1, x_2, x_3) &= f \big( [1]; x_1, g([1] + [2]; x_2, x_3) \big) - f \big( [2]; x_2, g([1]; x_1, x_3) \big), \\
(f \diamond_2 g)([2]; x_1, x_2, x_3) &= f \big( [2]; x_1, g([1]; x_2, x_3) \big) - f \big( [1]; x_2, g([1] + [2]; x_1, x_3) \big), \\
(f \diamond_2 g)([3]; x_1, x_2, x_3) &= f \big( [2]; x_1, g([2]; x_2, x_3) \big) - f \big( [2]; x_2, g([2]; x_1, x_3) \big).
\end{align*}
For $f = g = \pi$, the components $[1], [2], [3]$ of $\llbracket \pi, \pi \rrbracket_\pL = 2 (\pi \diamond_1 \pi - \pi \diamond_2 \pi)$ give back the identities (\ref{pre-l1}), (\ref{pre-l2}) and (\ref{pre-l3}).

\begin{thm}\label{thm-pl-gla}\cite[Theorem 3.7]{das-pre}
$(C^{\ast+1}_\pL (V,V), \llbracket \cdot, \cdot \rrbracket_\pL)$ is a graded Lie algebra. Two bilinear maps $\triangleleft, \triangleright$ on $V$ define a pre-Leibniz algebra if and only if the element $\pi \in C^2_\pL (V,V)$ given by
\[
\pi ([1]; x, y) = x \triangleleft y, \qquad \pi([2]; x, y) = x \triangleright y, \qquad \forall x, y \in V,
\]
satisfies $\llbracket \pi, \pi \rrbracket_\pL = 0$.
\end{thm}

The graded Lie algebra structure is obtained from Balavoine's graded Lie algebra \cite{bala-operad} of the doubled space $V \oplus V$ by an injective morphism of graded Lie algebras \cite{das-pre}.

\medskip

\noindent Cohomology. We now recall the cohomology of a pre-Leibniz algebra with coefficients in a representation \cite{das-pre}. Let $(\g, \triangleleft, \triangleright)$ be a pre-Leibniz algebra with Maurer-Cartan element $\pi$, and let $M$ be a representation. Define, $\forall x \in \g, ~ u \in M$,
\[
\begin{aligned}
\pi^L ([1]; x, u) &= x \triangleleft^L u, & \qquad \pi^L ([2]; x, u) &= x \triangleright^L u, \\
\pi^R ([1]; u, x) &= u \triangleleft^R x, & \qquad \pi^R ([2]; u, x) &= u \triangleright^R x.
\end{aligned}
\]
The $n$-th cochain group is $C^n_\pL (\g, M) = \mathrm{Hom}(\kk[C_n] \otimes \g^{\otimes n}, M)$, $n \geq 1$, and the coboundary $\delta_\pL : C^n_\pL (\g, M) \rightarrow C^{n+1}_\pL (\g, M)$ is given by \cite[(15)]{das-pre}
\begin{equation}\label{pl-diff}
\begin{split}
(\delta_\pL f)([r]; x_1, \ldots, x_{n+1}) = {}& \sum_{i=1}^n (-1)^{i+1} \pi^L \big( R_i [r]; x_i, f (S_i [r]; x_1, \ldots, \widehat{x_i}, \ldots, x_{n+1}) \big) \\
&{}+ (-1)^{n+1} \pi^R \big( R' [r]; f (S'[r]; x_1, \ldots, x_n), x_{n+1} \big) \\
&{}+ \sum_{1 \leq i < j \leq n+1} (-1)^i f \big( R_{ij}[r]; x_1, \ldots, \widehat{x_i}, \ldots, x_{j-1}, \\
&\hspace{5em} \pi (S_{ij}[r]; x_i, x_j), x_{j+1}, \ldots, x_{n+1} \big),
\end{split}
\end{equation}
for $[r] \in C_{n+1}$, where the pairs of maps $(R_i, S_i)$, $(R', S')$ and $(R_{ij}, S_{ij})$ are defined by the same rule as above (with the distinguished block containing the arguments of the inner operation) from the ordered partitions into blocks
\begin{gather*}
\big\{ [i] \big\}, ~ \big\{ [1], \ldots, \widehat{[i]}, \ldots, [n+1] \big\}; \qquad\qquad \big\{ [1], \ldots, [n] \big\}, ~ \big\{ [n+1] \big\}; \\*
\big\{ [1] \big\}, ~ \ldots, ~ \widehat{\big\{ [i] \big\}}, ~ \ldots, ~ \big\{ [j-1] \big\}, ~ \big\{ [i], [j] \big\}, ~ \big\{ [j+1] \big\}, ~ \ldots, ~ \big\{ [n+1] \big\},
\end{gather*}
respectively. (These are the blocks used in \cite[(15)]{das-pre}; we write them out explicitly, since the indices of $R$ and $S$ in that formula differ from the convention of \cite[(11)]{das-pre}.) For instance, for $n = 2$ and $(i,j) = (1,3)$, the blocks are $\{[2]\}, \{[1],[3]\}$, so that $R_{13}[2] = [1]$, $S_{13}[2] = [1]+[2]$, $R_{13}[1] = R_{13}[3] = [2]$, $S_{13}[1] = [1]$ and $S_{13}[3] = [2]$. Then $\delta_\pL \circ \delta_\pL = 0$ and the cohomology groups are denoted by $H^\ast_\pL (\g, M)$. For the adjoint representation,
\[
\delta_\pL f = (-1)^{n-1} \llbracket \pi, f \rrbracket_\pL, \qquad \forall f \in C^n_\pL (\g, \g).
\]
More generally, the coboundary with coefficients is induced by the semidirect product:

\begin{lemma}\label{lemma-delta-bracket}\cite[Proposition 4.3]{das-pre}
Let $\bar{\pi}$ be the Maurer-Cartan element of the semidirect product $\g \ltimes M$, and extend $f \in C^n_\pL (\g, M)$ to $\g \oplus M$ by zero whenever one of its arguments lies in $M$. Then $\llbracket \bar{\pi}, f \rrbracket_\pL$ belongs to $C^{n+1}_\pL (\g, M)$ and
\[
\delta_\pL f = (-1)^{n-1} \llbracket \bar{\pi}, f \rrbracket_\pL.
\]
\end{lemma}

This can also be checked directly by comparing (\ref{pl-bracket}) with (\ref{pl-diff}). For $n = 1$, writing $r_\triangleleft = 1$ and $r_\triangleright = 2$, the formula (\ref{pl-diff}) reads
\begin{equation}\label{pl-diff-1}
(\delta_\pL f)([r_\ast]; x, y) = x \ast f(y) + f(x) \ast y - f (x \ast y), \qquad \forall f \in \mathrm{Hom}(\g, M), ~ x, y \in \g, ~ \ast \in \{ \triangleleft, \triangleright \}.
\end{equation}

\begin{thm}\label{thm-phi}\cite[Theorem 4.5]{das-pre}
The maps
\[
\Phi_n : C^n_\pL (\g, M) \rightarrow C^n_\Leib (\g_\Tot, M_\Tot), \quad \Phi_n (f) (x_1, \ldots, x_n) = \sum_{r=1}^n f([r]; x_1, \ldots, x_n),
\]
define a morphism of cochain complexes. Hence there is a morphism $H^\ast_\pL (\g, M) \rightarrow H^\ast_\Leib (\g_\Tot, M_\Tot)$.
\end{thm}

\section{Weighted \texorpdfstring{$\mathcal{O}$}{O}-operators on pre-Leibniz algebras}\label{sec-3}

In this section, we introduce $\lambda$-weighted Rota-Baxter operators and $\lambda$-weighted $\mathcal{O}$-operators on pre-Leibniz algebras (see also \cite{sun-guo}) and study some basic properties.

\begin{defn}
Let $(\g, \triangleleft, \triangleright)$ be a pre-Leibniz algebra. A linear map $R : \g \rightarrow \g$ is said to be a $\lambda$-weighted Rota-Baxter operator on $\g$ if
\begin{align*}
R(x) \ast R(y) = R \big( R(x) \ast y + x \ast R(y) + \lambda \, x \ast y \big), \qquad \forall x, y \in \g, ~ \ast \in \{ \triangleleft, \triangleright \}.
\end{align*}
\end{defn}

\begin{exam}\label{exam-wrb}
(i) A Rota-Baxter operator on $\g$ is simply a $0$-weighted Rota-Baxter operator.

(ii) The identity map $\mathrm{id} : \g \rightarrow \g$ is a $(-1)$-weighted Rota-Baxter operator.

(iii) If $f : \g \rightarrow \g$ is a pre-Leibniz algebra morphism and $f^2 = f$ (idempotent condition), then $f$ is a $(-1)$-weighted Rota-Baxter operator.

(iv) If $R$ is a $\lambda$-weighted Rota-Baxter operator, then $\forall \nu \in \kk$, the map $\nu R$ is a $(\nu \lambda)$-weighted Rota-Baxter operator.

(v) If $R$ is a $\lambda$-weighted Rota-Baxter operator, then $- \lambda \, \mathrm{id} - R$ is a $\lambda$-weighted Rota-Baxter operator.
\end{exam}

\medskip

Let $\g$ and $\h$ be two pre-Leibniz algebras. We say that $\h$ is a pre-Leibniz $\g$-representation if there are bilinear maps
\[
\triangleleft^L, \triangleright^L : \g \otimes \h \rightarrow \h, \qquad \triangleleft^R, \triangleright^R : \h \otimes \g \rightarrow \h
\]
that make $(\h, \triangleleft^L, \triangleright^L, \triangleleft^R, \triangleright^R)$ into a representation of the pre-Leibniz algebra $\g$ (Definition \ref{defn-rep}) satisfying additionally the following nine identities:
\begin{align}
\shortintertext{$\forall x \in \g,~ y, z \in \h$:}
x \triangleleft^L (y \triangleleft_\h z + y \triangleright_\h z) &= (x \triangleleft^L y) \triangleleft_\h z + y \triangleright_\h (x \triangleleft^L z) \label{comp-1} \\
x \triangleright^L (y \triangleleft_\h z) &= (x \triangleright^L y) \triangleleft_\h z + y \triangleleft_\h (x \triangleleft^L z + x \triangleright^L z) \label{comp-2} \\
x \triangleright^L (y \triangleright_\h z) &= (x \triangleleft^L y + x \triangleright^L y) \triangleright_\h z + y \triangleright_\h (x \triangleright^L z) \label{comp-3} \\
\shortintertext{$\forall y \in \g,~ x, z \in \h$:}
x \triangleleft_\h (y \triangleleft^L z + y \triangleright^L z) &= (x \triangleleft^R y) \triangleleft_\h z + y \triangleright^L (x \triangleleft_\h z) \label{comp-4} \\
x \triangleright_\h (y \triangleleft^L z) &= (x \triangleright^R y) \triangleleft_\h z + y \triangleleft^L (x \triangleleft_\h z + x \triangleright_\h z) \label{comp-5} \\
x \triangleright_\h (y \triangleright^L z) &= (x \triangleleft^R y + x \triangleright^R y) \triangleright_\h z + y \triangleright^L (x \triangleright_\h z) \label{comp-6} \\
\shortintertext{$\forall z \in \g,~ x, y \in \h$:}
x \triangleleft_\h (y \triangleleft^R z + y \triangleright^R z) &= (x \triangleleft_\h y) \triangleleft^R z + y \triangleright_\h (x \triangleleft^R z) \label{comp-7} \\
x \triangleright_\h (y \triangleleft^R z) &= (x \triangleright_\h y) \triangleleft^R z + y \triangleleft_\h (x \triangleleft^R z + x \triangleright^R z) \label{comp-8} \\
x \triangleright_\h (y \triangleright^R z) &= (x \triangleleft_\h y + x \triangleright_\h y) \triangleright^R z + y \triangleright_\h (x \triangleright^R z) \label{comp-9}
\end{align}
Note that, for any pre-Leibniz algebra $\g$, the adjoint representation is a pre-Leibniz $\g$-representation. Thus a pre-Leibniz $\g$-representation is a pre-Leibniz algebra $\h$ together with a representation of $\g$ on it, compatible with the operations of $\h$; when $\h$ carries the zero operations, it is simply a representation of $\g$ (Definition \ref{defn-rep}). Equivalently, $\h$ is a pre-Leibniz $\g$-representation if and only if the direct sum $\g \oplus \h$ is a pre-Leibniz algebra with the operations
\begin{equation}\label{wsd-1}
(x, u) \ast (y, v) = \big( x \ast y, ~ x \ast^L v + u \ast^R y + u \ast_\h v \big), \qquad \forall x, y \in \g, ~ u, v \in \h, ~ \ast \in \{ \triangleleft, \triangleright \}.
\end{equation}
Indeed, the identities (\ref{pre-l1})-(\ref{pre-l3}) for $\g \oplus \h$ split according to the number of arguments taken in $\h$: with no argument in $\h$, they are the identities of $\g$; with one argument in $\h$, they are the axioms of a representation of $\g$; with two arguments in $\h$, they are the nine identities (\ref{comp-1})-(\ref{comp-9}); with three arguments in $\h$, they are the identities of $\h$.

\begin{defn}\label{defn-wrrb}
Let $\g$ and $\h$ be two pre-Leibniz algebras and $\h$ be a pre-Leibniz $\g$-representation. A linear map $T : \h \rightarrow \g$ is said to be a $\lambda$-weighted $\mathcal{O}$-operator (on $\h$ over the pre-Leibniz algebra $\g$) if
\begin{equation}\label{lw}
T(u) \ast T(v) = T \big( T(u) \ast^L v + u \ast^R T(v) + \lambda \, u \ast_\h v \big), \qquad \forall u, v \in \h, ~ \ast \in \{ \triangleleft, \triangleright \}.
\end{equation}
\end{defn}

It follows that a $\lambda$-weighted Rota-Baxter operator on a pre-Leibniz algebra $\g$ is a $\lambda$-weighted $\mathcal{O}$-operator on $\g$ over itself.

In addition to Example \ref{exam-wrb}, we have the following examples.

\begin{exam}\label{exam-wrrb}
(i) Let $\g, \h$ be two pre-Leibniz algebras and $\h$ be a pre-Leibniz $\g$-representation. A linear map $D : \g \rightarrow \h$ is said to be a $\lambda$-weighted crossed homomorphism if
\[
D (x \ast y) = x \ast^L D(y) + D(x) \ast^R y + \lambda \, D(x) \ast_\h D(y), \qquad \forall x, y \in \g, ~ \ast \in \{ \triangleleft, \triangleright \}.
\]
If $D$ is an invertible $\lambda$-weighted crossed homomorphism, then $D^{-1} : \h \rightarrow \g$ is a $\lambda$-weighted $\mathcal{O}$-operator.

(ii) Let $\g$ be a pre-Leibniz algebra. An ideal of $\g$ is a subspace $\h$ such that $\g \ast \h \subset \h$ and $\h \ast \g \subset \h$, $\forall \ast \in \{ \triangleleft, \triangleright \}$. Any ideal $\h$ of $\g$ can be regarded as a pre-Leibniz $\g$-representation, where the operations of $\h$ and the four actions are the restrictions of $\triangleleft$ and $\triangleright$. Then the inclusion map $i : \h \rightarrow \g$ is a $(-1)$-weighted $\mathcal{O}$-operator.
\end{exam}

To characterize $\lambda$-weighted $\mathcal{O}$-operators, we first need a weighted version of the semidirect product.

\begin{prop}\label{prop-wsd}
(Weighted semidirect product) Let $\g$ and $\h$ be two pre-Leibniz algebras and $\h$ be a pre-Leibniz $\g$-representation. Then $\forall \lambda \in \kk$, the direct sum $\g \oplus \h$ carries a new pre-Leibniz algebra structure with operations
\[
(x, u) \ast_\lambda (y, v) = \big( x \ast y, ~ x \ast^L v + u \ast^R y + \lambda \, u \ast_\h v \big), \qquad \forall (x,u), (y,v) \in \g \oplus \h, ~ \ast \in \{ \triangleleft, \triangleright \}.
\]
This is called the $\lambda$-weighted semidirect product, often denoted by $\g \ltimes_\lambda \h$.
\end{prop}

\begin{proof}
For $\lambda = 0$, this is the semidirect product $\g \ltimes \h$ of Definition \ref{defn-rep}. For $\lambda \neq 0$, the linear isomorphism $\psi_\lambda (x, u) = (x, \lambda u)$ satisfies
\[
\psi_\lambda \big( (x,u) \ast_\lambda (y,v) \big) = \big( x \ast y, ~ \lambda \, x \ast^L v + \lambda \, u \ast^R y + \lambda^2 u \ast_\h v \big) = \psi_\lambda (x,u) \ast \psi_\lambda (y,v),
\]
where the right hand side is computed with the operations (\ref{wsd-1}). Since $\g \oplus \h$ with the operations (\ref{wsd-1}) is a pre-Leibniz algebra, so is $\g \ltimes_\lambda \h$.
\end{proof}

We have the following characterization of $\lambda$-weighted $\mathcal{O}$-operators.

\begin{prop}\label{prop-graph}
A linear map $T : \h \rightarrow \g$ is a $\lambda$-weighted $\mathcal{O}$-operator if and only if the graph $Gr (T) = \{ (T(u), u) ~|~ u \in \h \}$ is a subalgebra of the $\lambda$-weighted semidirect product $\g \ltimes_\lambda \h$.
\end{prop}

\begin{proof}
$\forall u, v \in \h$ and $\ast \in \{ \triangleleft, \triangleright \}$, we have
\begin{align*}
(T(u), u) \ast_\lambda (T(v), v) = \big( T(u) \ast T(v), ~ T(u) \ast^L v + u \ast^R T(v) + \lambda \, u \ast_\h v \big).
\end{align*}
This is in $Gr (T)$ if and only if (\ref{lw}) holds. Hence the result follows.
\end{proof}

\begin{remark}\label{rmk-induced}
Let $T : \h \rightarrow \g$ be a $\lambda$-weighted $\mathcal{O}$-operator. Since $Gr (T)$ is isomorphic to $\h$ as a vector space, we get that $\h$ inherits a new pre-Leibniz algebra structure with the operations
\begin{equation}\label{induced-pl}
u \ast_T v := T(u) \ast^L v + u \ast^R T(v) + \lambda \, u \ast_\h v, \qquad \forall u, v \in \h, ~ \ast \in \{ \triangleleft, \triangleright \}.
\end{equation}
In other words, $(\h, \triangleleft_T, \triangleright_T)$ is a pre-Leibniz algebra, denoted by $\h_T$ (called the induced pre-Leibniz algebra). Moreover, $T : \h_T \rightarrow \g$ is a morphism of pre-Leibniz algebras, since $T (u \ast_T v) = T(u) \ast T(v)$ by (\ref{lw}).
\end{remark}

Let $T, T' : \h \rightarrow \g$ be two $\lambda$-weighted $\mathcal{O}$-operators. A morphism from $T$ to $T'$ consists of a pair $(\phi, \psi)$ of pre-Leibniz algebra morphisms $\phi : \g \rightarrow \g$ and $\psi : \h \rightarrow \h$ satisfying
\[
\phi \circ T = T' \circ \psi, \qquad \psi (x \ast^L u) = \phi (x) \ast^L \psi (u), \qquad \psi (u \ast^R x) = \psi (u) \ast^R \phi (x),
\]
$\forall x \in \g$, $u \in \h$ and $\ast \in \{ \triangleleft, \triangleright \}$.

\begin{prop}\label{prop-morph}
Let $(\phi, \psi)$ be a morphism of $\lambda$-weighted $\mathcal{O}$-operators from $T$ to $T'$. Then $\psi : \h \rightarrow \h$ is a morphism of induced pre-Leibniz algebras from $\h_T$ to $\h_{T'}$.
\end{prop}

\begin{proof}
$\forall u, v \in \h$ and $\ast \in \{ \triangleleft, \triangleright \}$, we have
\begin{align*}
\psi (u \ast_T v) &= \psi \big( T(u) \ast^L v + u \ast^R T(v) + \lambda \, u \ast_\h v \big) \\
&= \phi (T(u)) \ast^L \psi (v) + \psi (u) \ast^R \phi (T(v)) + \lambda \, \psi (u) \ast_\h \psi (v) \\
&= T'(\psi (u)) \ast^L \psi (v) + \psi (u) \ast^R T'(\psi (v)) + \lambda \, \psi (u) \ast_\h \psi (v) = \psi (u) \ast_{T'} \psi (v).
\end{align*}
This shows that $\psi : \h_T \rightarrow \h_{T'}$ is a morphism of pre-Leibniz algebras.
\end{proof}

\section{Maurer-Cartan characterization and cohomology of weighted \texorpdfstring{$\mathcal{O}$}{O}-operators}\label{sec-4}

Let $\g$ and $\h$ be two pre-Leibniz algebras and $\h$ be a pre-Leibniz $\g$-representation. With this data, we first construct a differential graded Lie algebra whose Maurer-Cartan elements are precisely $\lambda$-weighted $\mathcal{O}$-operators. This characterization allows us to define a cohomology associated to a $\lambda$-weighted $\mathcal{O}$-operator. Later, we interpret this as the cohomology of the corresponding induced pre-Leibniz algebra with coefficients in a suitable representation.

\medskip

Let $\pi_\g$ and $\pi_\h$ denote the Maurer-Cartan elements corresponding to the pre-Leibniz algebras $\g$ and $\h$, respectively. Take $V = \g \oplus \h$ and consider the graded Lie algebra $(C^{\ast+1}_\pL (V, V), \llbracket \cdot, \cdot \rrbracket_\pL)$ of Theorem \ref{thm-pl-gla}. Then the element $\theta \in C^2_\pL (V,V)$ defined by
\[
\theta \big( [r_\ast]; (x,u), (y,v) \big) = \big( x \ast y, ~ x \ast^L v + u \ast^R y \big), \qquad \forall (x,u), (y,v) \in V, ~ \ast \in \{ \triangleleft, \triangleright \}
\]
(recall $r_\triangleleft = 1$ and $r_\triangleright = 2$) is the Maurer-Cartan element of the semidirect product $\g \ltimes \h$, so it satisfies $\llbracket \theta, \theta \rrbracket_\pL = 0$. Hence $\theta$ induces a degree $1$ differential $d_\theta = \llbracket \theta, - \rrbracket_\pL$ on $C^{\ast+1}_\pL (V, V)$. Moreover, the graded subspace
\[
\mathcal{A} = \bigoplus_{n \geq 1} \mathrm{Hom}(\kk[C_n] \otimes \h^{\otimes n}, \g) \subset C^{\ast}_\pL (V,V)
\]
(an element of $\mathrm{Hom}(\kk[C_n] \otimes \h^{\otimes n}, \g)$ being extended to $V$ by zero whenever one of its arguments lies in $\g$) is an abelian Lie subalgebra. Hence, by the derived bracket construction \cite{voro}, the graded space $\mathcal{A}$ carries a graded Lie bracket given by
\begin{equation}\label{derived}
\llbracket P, Q \rrbracket := (-1)^m \llbracket d_\theta (P), Q \rrbracket_\pL = (-1)^m \llbracket \llbracket \theta, P \rrbracket_\pL, Q \rrbracket_\pL,
\end{equation}
$\forall P \in \mathrm{Hom}(\kk[C_m] \otimes \h^{\otimes m}, \g)$ and $Q \in \mathrm{Hom}(\kk[C_n] \otimes \h^{\otimes n}, \g)$. Here $\mathrm{Hom}(\kk[C_n] \otimes \h^{\otimes n}, \g)$ has degree $n$ in $\mathcal{A}$ (while it has degree $n-1$ in $C^{\ast+1}_\pL (V,V)$); in particular, a linear map $T : \h \rightarrow \g$ has degree $1$. The explicit description of the bracket is as follows (compare with \cite[(10)]{das-wrb}): $\forall [r] \in C_{m+n}$ and $u_1, \ldots, u_{m+n} \in \h$,
\begin{align}\label{bracket-explicit}
&\llbracket P, Q \rrbracket ([r]; u_1, \ldots, u_{m+n}) \\
&= \sum_{i=1}^m \sum_{\sigma \in Sh(i-1, n)} (-1)^{(i-1)n} (-1)^\sigma P \big( u_{\sigma(1)}, \ldots, u_{\sigma(i-1)}, \nonumber \\
&\hspace{9em} Q (u_{\sigma(i)}, \ldots, u_{\sigma(i+n-1)}) \ast^L u_{i+n}, u_{i+n+1}, \ldots, u_{m+n} \big) \nonumber \\
&+ \sum_{i=1}^m \sum_{\sigma \in Sh(i-1, 1, n-1)} (-1)^{(i-1)n + n - 1} (-1)^\sigma P \big( u_{\sigma(1)}, \ldots, u_{\sigma(i-1)}, \nonumber \\
&\hspace{9em} u_{\sigma(i)} \ast^R Q (u_{\sigma(i+1)}, \ldots, u_{\sigma(i+n-1)}, u_{i+n}), u_{i+n+1}, \ldots, u_{m+n} \big) \nonumber \\
&+ \sum_{\sigma \in Sh(m, n-1)} (-1)^{mn} (-1)^\sigma P (u_{\sigma(1)}, \ldots, u_{\sigma(m)}) \ast Q (u_{\sigma(m+1)}, \ldots, u_{\sigma(m+n-1)}, u_{m+n}) \nonumber \\
&- (-1)^{mn} \Big\{ \sum_{i=1}^n \sum_{\sigma \in Sh(i-1, m)} (-1)^{(i-1)m} (-1)^\sigma Q \big( u_{\sigma(1)}, \ldots, u_{\sigma(i-1)}, \nonumber \\
&\hspace{9em} P (u_{\sigma(i)}, \ldots, u_{\sigma(i+m-1)}) \ast^L u_{i+m}, u_{i+m+1}, \ldots, u_{m+n} \big) \nonumber \\
&+ \sum_{i=1}^n \sum_{\sigma \in Sh(i-1, 1, m-1)} (-1)^{(i-1)m + m - 1} (-1)^\sigma Q \big( u_{\sigma(1)}, \ldots, u_{\sigma(i-1)}, \nonumber \\
&\hspace{9em} u_{\sigma(i)} \ast^R P (u_{\sigma(i+1)}, \ldots, u_{\sigma(i+m-1)}, u_{i+m}), u_{i+m+1}, \ldots, u_{m+n} \big) \nonumber \\
&+ \sum_{\sigma \in Sh(n, m-1)} (-1)^{mn} (-1)^\sigma Q (u_{\sigma(1)}, \ldots, u_{\sigma(n)}) \ast P (u_{\sigma(n+1)}, \ldots, u_{\sigma(m+n-1)}, u_{m+n}) \Big\}. \nonumber
\end{align}
Here $Sh(p, q)$ (resp. $Sh(p, q, s)$) denotes the set of $(p,q)$-shuffles (resp. $(p,q,s)$-shuffles), and the labels are determined by $[r]$, according to the rule of the maps $R$ and $S$: each of the operations $P$, $Q$ and $\ast$ occurring in a term receives the label $[k]$ if the variable $u_r$ appears in its $k$-th argument, and the sum of all its labels otherwise. In particular, $\ast$ stands for $\triangleleft$ (label $[1]$), for $\triangleright$ (label $[2]$), or for $\triangleleft + \triangleright$ (sum of the labels); similarly $\ast^L \in \{ \triangleleft^L, \triangleright^L, \triangleleft^L + \triangleright^L \}$ and $\ast^R \in \{ \triangleleft^R, \triangleright^R, \triangleleft^R + \triangleright^R \}$.

On the other hand, $\forall \lambda \in \kk$, the element $\theta' = - \lambda \pi_\h \in C^2_\pL (V,V)$ defined by
\[
\theta' \big( [r_\ast]; (x,u), (y,v) \big) = \big( 0, ~ - \lambda \, u \ast_\h v \big), \qquad \forall (x,u), (y,v) \in V, ~ \ast \in \{ \triangleleft, \triangleright \}
\]
satisfies $\llbracket \theta', \theta' \rrbracket_\pL = 0$. Hence the element $\theta'$ also induces a differential $d_{\theta'} = \llbracket \theta', - \rrbracket_\pL$ on $C^{\ast+1}_\pL (V,V)$. It is easy to see that the graded subspace $\mathcal{A}$ is closed under the differential $d_{\theta'}$. Let us denote the restriction of $d_{\theta'}$ to this subspace by $d$. This is explicitly given by
\begin{equation}\label{d-formula}
\begin{split}
(dP)([r]; u_1, \ldots, u_{n+1}) = {}& (-1)^n \sum_{1 \leq i < j \leq n+1} (-1)^i P \big( R_{ij}[r]; u_1, \ldots, \widehat{u_i}, \ldots, u_{j-1}, \\
&\hspace{5em} \lambda \, \pi_\h (S_{ij}[r]; u_i, u_j), u_{j+1}, \ldots, u_{n+1} \big),
\end{split}
\end{equation}
with the notations of (\ref{pl-diff}).

Finally, we observe that the elements $\theta$ and $\theta'$ additionally satisfy the compatibility $\llbracket \theta, \theta' \rrbracket_\pL = 0$. Indeed, $\theta + \pi_\h$ is the Maurer-Cartan element of $\g \ltimes_1 \h$, so $0 = \llbracket \theta + \pi_\h, \theta + \pi_\h \rrbracket_\pL = 2 \llbracket \theta, \pi_\h \rrbracket_\pL$. This implies that
\begin{align*}
d \llbracket P, Q \rrbracket &= (-1)^m \llbracket \theta', \llbracket \llbracket \theta, P \rrbracket_\pL, Q \rrbracket_\pL \rrbracket_\pL \\
&= (-1)^m \llbracket \llbracket \theta', \llbracket \theta, P \rrbracket_\pL \rrbracket_\pL, Q \rrbracket_\pL + (-1)^m (-1)^m \llbracket \llbracket \theta, P \rrbracket_\pL, \llbracket \theta', Q \rrbracket_\pL \rrbracket_\pL \\
&= (-1)^{m+1} \llbracket \llbracket \theta, \llbracket \theta', P \rrbracket_\pL \rrbracket_\pL, Q \rrbracket_\pL + \llbracket \llbracket \theta, P \rrbracket_\pL, dQ \rrbracket_\pL \\
&= (-1)^{m+1} \llbracket \llbracket \theta, dP \rrbracket_\pL, Q \rrbracket_\pL + \llbracket \llbracket \theta, P \rrbracket_\pL, dQ \rrbracket_\pL = \llbracket dP, Q \rrbracket + (-1)^m \llbracket P, dQ \rrbracket.
\end{align*}
Thus, we get the following.

\begin{thm}\label{mc-thm}
Let $\g, \h$ be two pre-Leibniz algebras and $\h$ be a pre-Leibniz $\g$-representation. Then the triple $(\mathcal{A}, \llbracket \cdot, \cdot \rrbracket, d)$ is a differential graded Lie algebra. Moreover, a linear map $T : \h \rightarrow \g$ is a $\lambda$-weighted $\mathcal{O}$-operator if and only if $T$ is a Maurer-Cartan element in the differential graded Lie algebra $(\mathcal{A}, \llbracket \cdot, \cdot \rrbracket, d)$.
\end{thm}

\begin{proof}
The first part follows from previous discussions. For the second part, we observe that for any linear map $T : \h \rightarrow \g$, we have, $\forall u, v \in \h$, $\ast \in \{ \triangleleft, \triangleright \}$,
\begin{align*}
\llbracket T, T \rrbracket ([r_\ast]; u, v) &= 2 \big( T \big( T(u) \ast^L v + u \ast^R T(v) \big) - T(u) \ast T(v) \big), \\
(dT)([r_\ast]; u, v) &= \lambda \, T (u \ast_\h v),
\end{align*}
by (\ref{bracket-explicit}) with $m = n = 1$ and (\ref{d-formula}) with $(i,j) = (1,2)$. Therefore,
\[
\big( dT + \tfrac{1}{2} \llbracket T, T \rrbracket \big)([r_\ast]; u, v) = T \big( T(u) \ast^L v + u \ast^R T(v) + \lambda \, u \ast_\h v \big) - T(u) \ast T(v).
\]
This shows that $T$ is a Maurer-Cartan element (i.e. $dT + \frac{1}{2} \llbracket T, T \rrbracket = 0$) if and only if $T$ is a $\lambda$-weighted $\mathcal{O}$-operator.
\end{proof}

Let $T : \h \rightarrow \g$ be a $\lambda$-weighted $\mathcal{O}$-operator. That is, $T$ can be equivalently seen as a Maurer-Cartan element in the differential graded Lie algebra constructed in Theorem \ref{mc-thm}. Therefore, $T$ induces a differential
\[
d_T = d + \llbracket T, - \rrbracket : \mathrm{Hom}(\kk[C_n] \otimes \h^{\otimes n}, \g) \rightarrow \mathrm{Hom}(\kk[C_{n+1}] \otimes \h^{\otimes n+1}, \g), \qquad \forall n \geq 1.
\]
We set, $\forall u, v \in \h$, $x \in \g$ and $\ast \in \{ \triangleleft, \triangleright \}$,
\begin{gather*}
\pi^L_T ([r_\ast]; u, x) = T(u) \ast x - T(u \ast^R x), \qquad \pi^R_T ([r_\ast]; x, u) = x \ast T(u) - T(x \ast^L u), \\
\pi_T ([r_\ast]; u, v) = u \ast_T v = T(u) \ast^L v + u \ast^R T(v) + \lambda \, u \ast_\h v.
\end{gather*}
Explicitly, $\forall f \in \mathrm{Hom}(\kk[C_n] \otimes \h^{\otimes n}, \g)$, $[r] \in C_{n+1}$ and $u_1, \ldots, u_{n+1} \in \h$, we have
\begin{equation}\label{dT-explicit}
\begin{split}
(d_T f)([r]; u_1, \ldots, u_{n+1}) = {}& (-1)^n \sum_{i=1}^n (-1)^{i+1} \pi^L_T \big( R_i[r]; u_i, f (S_i[r]; u_1, \ldots, \widehat{u_i}, \ldots, u_{n+1}) \big) \\
&{}- \pi^R_T \big( R'[r]; f (S'[r]; u_1, \ldots, u_n), u_{n+1} \big) \\
&{}+ (-1)^n \sum_{1 \leq i < j \leq n+1} (-1)^i f \big( R_{ij}[r]; u_1, \ldots, \widehat{u_i}, \ldots, u_{j-1}, \\
&\hspace{5em} \pi_T (S_{ij}[r]; u_i, u_j), u_{j+1}, \ldots, u_{n+1} \big),
\end{split}
\end{equation}
where the pairs $(R_i, S_i)$, $(R', S')$ and $(R_{ij}, S_{ij})$ are as in (\ref{pl-diff}). This follows from (\ref{bracket-explicit}) with $m = 1$, $P = T$ and $Q = f$, together with (\ref{d-formula}).

In the following, we show that the differential (\ref{dT-explicit}) can be seen as the differential of the pre-Leibniz algebra cohomology of $\h_T$ (see Remark \ref{rmk-induced}) with coefficients in a suitable representation. We first observe the following.

\begin{prop}\label{prop-rep-T}
Let $T : \h \rightarrow \g$ be a $\lambda$-weighted $\mathcal{O}$-operator. Then $\g$ is a representation of the induced pre-Leibniz algebra $\h_T$ with the actions given by
\[
u \ast^L_T x = T(u) \ast x - T (u \ast^R x), \qquad x \ast^R_T u = x \ast T(u) - T (x \ast^L u), \qquad \forall u \in \h, ~ x \in \g, ~ \ast \in \{ \triangleleft, \triangleright \}.
\]
\end{prop}

\begin{proof}
We have to verify the nine identities of a representation (Definition \ref{defn-rep}) for the pre-Leibniz algebra $\h_T$ acting on $\g$. Let $u, v \in \h$ and $x \in \g$. Recall that $T (u \ast_T v) = T(u) \ast T(v)$ by (\ref{lw}). We first verify the identity (\ref{pre-l1}) with $x$ in the third position. Using the definition of the actions, we get
\begin{align*}
& u \triangleleft^L_T \big( v \triangleleft^L_T x + v \triangleright^L_T x \big) - (u \triangleleft_T v) \triangleleft^L_T x - v \triangleright^L_T \big( u \triangleleft^L_T x \big) \\
&= \cancel{T(u) \triangleleft \big( T(v) \triangleleft x + T(v) \triangleright x \big)} - T(u) \triangleleft T(w) \\
&\quad - T \big( u \triangleleft^R ( T(v) \triangleleft x + T(v) \triangleright x ) \big) + T \big( u \triangleleft^R T(w) \big) \\
&\quad - \cancel{\big( T(u) \triangleleft T(v) \big) \triangleleft x} + T \big( ( T(u) \triangleleft^L v + u \triangleleft^R T(v) + \lambda \, u \triangleleft_\h v ) \triangleleft^R x \big) \\
&\quad - \cancel{T(v) \triangleright \big( T(u) \triangleleft x \big)} + T(v) \triangleright T \big( u \triangleleft^R x \big) + T \big( v \triangleright^R ( T(u) \triangleleft x ) \big) - T \big( v \triangleright^R T ( u \triangleleft^R x ) \big),
\end{align*}
where $w = v \triangleleft^R x + v \triangleright^R x$, and the crossed terms cancel by the identity (\ref{pre-l1}) in $\g$. Expanding $T(u) \triangleleft T(w)$ and $T(v) \triangleright T( u \triangleleft^R x )$ by (\ref{lw}), the above expression is equal to
\begin{align*}
& - T \big( T(u) \triangleleft^L w + \cancel{u \triangleleft^R T(w)} + \lambda \, u \triangleleft_\h w \big) \\
&\quad - T \big( u \triangleleft^R ( T(v) \triangleleft x + T(v) \triangleright x ) \big) + \cancel{T \big( u \triangleleft^R T(w) \big)} \\
&\quad + T \big( ( T(u) \triangleleft^L v + u \triangleleft^R T(v) + \lambda \, u \triangleleft_\h v ) \triangleleft^R x \big) \\
&\quad + T \big( T(v) \triangleright^L ( u \triangleleft^R x ) + \cancel{v \triangleright^R T ( u \triangleleft^R x )} + \lambda \, v \triangleright_\h ( u \triangleleft^R x ) \big) \\
&\quad + T \big( v \triangleright^R ( T(u) \triangleleft x ) \big) - \cancel{T \big( v \triangleright^R T ( u \triangleleft^R x ) \big)} \\
&= - T \Big( T(u) \triangleleft^L \big( v \triangleleft^R x + v \triangleright^R x \big) - \big( T(u) \triangleleft^L v \big) \triangleleft^R x - v \triangleright^R \big( T(u) \triangleleft x \big) \Big) \\
&\quad - T \Big( u \triangleleft^R \big( T(v) \triangleleft x + T(v) \triangleright x \big) - \big( u \triangleleft^R T(v) \big) \triangleleft^R x - T(v) \triangleright^L \big( u \triangleleft^R x \big) \Big) \\
&\quad - \lambda \, T \Big( u \triangleleft_\h \big( v \triangleleft^R x + v \triangleright^R x \big) - \big( u \triangleleft_\h v \big) \triangleleft^R x - v \triangleright_\h \big( u \triangleleft^R x \big) \Big) \\
&= 0.
\end{align*}
Indeed, the first two brackets vanish by the identity (\ref{pre-l1}) for the representation $\h$ of $\g$ (with arguments $(T(u), v, x)$ and $(u, T(v), x)$), and the third one vanishes by the identity (\ref{comp-7}). Similarly, for the identity (\ref{pre-l2}) with $x$ in the third position, we get
\begin{align*}
& u \triangleright^L_T \big( v \triangleleft^L_T x \big) - (u \triangleright_T v) \triangleleft^L_T x - v \triangleleft^L_T \big( u \triangleleft^L_T x + u \triangleright^L_T x \big) \\
&= \cancel{T(u) \triangleright \big( T(v) \triangleleft x \big)} - T(u) \triangleright T \big( v \triangleleft^R x \big) - T \big( u \triangleright^R ( T(v) \triangleleft x ) \big) + T \big( u \triangleright^R T ( v \triangleleft^R x ) \big) \\
&\quad - \cancel{\big( T(u) \triangleright T(v) \big) \triangleleft x} + T \big( ( T(u) \triangleright^L v + u \triangleright^R T(v) + \lambda \, u \triangleright_\h v ) \triangleleft^R x \big) \\
&\quad - \cancel{T(v) \triangleleft \big( T(u) \triangleleft x + T(u) \triangleright x \big)} + T(v) \triangleleft T(w') \\
&\quad + T \big( v \triangleleft^R ( T(u) \triangleleft x + T(u) \triangleright x ) \big) - T \big( v \triangleleft^R T(w') \big),
\end{align*}
where $w' = u \triangleleft^R x + u \triangleright^R x$, and the crossed terms cancel by the identity (\ref{pre-l2}) in $\g$. Expanding $T(u) \triangleright T(v \triangleleft^R x)$ and $T(v) \triangleleft T(w')$ by (\ref{lw}), the above expression is equal to
\begin{align*}
& - T \big( T(u) \triangleright^L ( v \triangleleft^R x ) + \cancel{u \triangleright^R T ( v \triangleleft^R x )} + \lambda \, u \triangleright_\h ( v \triangleleft^R x ) \big) \\
&\quad - T \big( u \triangleright^R ( T(v) \triangleleft x ) \big) + \cancel{T \big( u \triangleright^R T ( v \triangleleft^R x ) \big)} \\
&\quad + T \big( ( T(u) \triangleright^L v + u \triangleright^R T(v) + \lambda \, u \triangleright_\h v ) \triangleleft^R x \big) \\
&\quad + T \big( T(v) \triangleleft^L w' + \cancel{v \triangleleft^R T(w')} + \lambda \, v \triangleleft_\h w' \big) \\
&\quad + T \big( v \triangleleft^R ( T(u) \triangleleft x + T(u) \triangleright x ) \big) - \cancel{T \big( v \triangleleft^R T(w') \big)} \\
&= - T \Big( T(u) \triangleright^L \big( v \triangleleft^R x \big) - \big( T(u) \triangleright^L v \big) \triangleleft^R x - v \triangleleft^R \big( T(u) \triangleleft x + T(u) \triangleright x \big) \Big) \\
&\quad - T \Big( u \triangleright^R \big( T(v) \triangleleft x \big) - \big( u \triangleright^R T(v) \big) \triangleleft^R x - T(v) \triangleleft^L \big( u \triangleleft^R x + u \triangleright^R x \big) \Big) \\
&\quad - \lambda \, T \Big( u \triangleright_\h \big( v \triangleleft^R x \big) - \big( u \triangleright_\h v \big) \triangleleft^R x - v \triangleleft_\h \big( u \triangleleft^R x + u \triangleright^R x \big) \Big) \\
&= 0,
\end{align*}
by the identity (\ref{pre-l2}) for the representation $\h$ and the identity (\ref{comp-8}). For the identity (\ref{pre-l3}) with $x$ in the third position, we get in the same way
\begin{align*}
& u \triangleright^L_T \big( v \triangleright^L_T x \big) - \big( u \triangleleft_T v + u \triangleright_T v \big) \triangleright^L_T x - v \triangleright^L_T \big( u \triangleright^L_T x \big) \\
&= \cancel{T(u) \triangleright \big( T(v) \triangleright x \big)} - T(u) \triangleright T \big( v \triangleright^R x \big) - T \big( u \triangleright^R ( T(v) \triangleright x ) \big) + T \big( u \triangleright^R T ( v \triangleright^R x ) \big) \\
&\quad - \cancel{\big( T(u) \triangleleft T(v) + T(u) \triangleright T(v) \big) \triangleright x} + T \big( ( u \triangleleft_T v + u \triangleright_T v ) \triangleright^R x \big) \\
&\quad - \cancel{T(v) \triangleright \big( T(u) \triangleright x \big)} + T(v) \triangleright T \big( u \triangleright^R x \big) + T \big( v \triangleright^R ( T(u) \triangleright x ) \big) - T \big( v \triangleright^R T ( u \triangleright^R x ) \big),
\end{align*}
where the crossed terms cancel by the identity (\ref{pre-l3}) in $\g$. Expanding $T(u) \triangleright T(v \triangleright^R x)$ and $T(v) \triangleright T(u \triangleright^R x)$ by (\ref{lw}), the above expression is equal to
\begin{align*}
& - T \big( T(u) \triangleright^L ( v \triangleright^R x ) + \cancel{u \triangleright^R T ( v \triangleright^R x )} + \lambda \, u \triangleright_\h ( v \triangleright^R x ) \big) \\
&\quad - T \big( u \triangleright^R ( T(v) \triangleright x ) \big) + \cancel{T \big( u \triangleright^R T ( v \triangleright^R x ) \big)} \\
&\quad + T \big( ( T(u) \triangleleft^L v + T(u) \triangleright^L v + u \triangleleft^R T(v) + u \triangleright^R T(v) + \lambda \, u \triangleleft_\h v + \lambda \, u \triangleright_\h v ) \triangleright^R x \big) \\
&\quad + T \big( T(v) \triangleright^L ( u \triangleright^R x ) + \cancel{v \triangleright^R T ( u \triangleright^R x )} + \lambda \, v \triangleright_\h ( u \triangleright^R x ) \big) \\
&\quad + T \big( v \triangleright^R ( T(u) \triangleright x ) \big) - \cancel{T \big( v \triangleright^R T ( u \triangleright^R x ) \big)} \\
&= - T \Big( T(u) \triangleright^L \big( v \triangleright^R x \big) - \big( T(u) \triangleleft^L v + T(u) \triangleright^L v \big) \triangleright^R x - v \triangleright^R \big( T(u) \triangleright x \big) \Big) \\
&\quad - T \Big( u \triangleright^R \big( T(v) \triangleright x \big) - \big( u \triangleleft^R T(v) + u \triangleright^R T(v) \big) \triangleright^R x - T(v) \triangleright^L \big( u \triangleright^R x \big) \Big) \\
&\quad - \lambda \, T \Big( u \triangleright_\h \big( v \triangleright^R x \big) - \big( u \triangleleft_\h v + u \triangleright_\h v \big) \triangleright^R x - v \triangleright_\h \big( u \triangleright^R x \big) \Big) \\
&= 0,
\end{align*}
by the identity (\ref{pre-l3}) for the representation $\h$ and the identity (\ref{comp-9}).
The six identities in which $x$ is the first or the second argument are verified in the same way.
%The six identities in which $x$ is the first or the second argument are verified by the same computation: the terms without $T$ cancel by (\ref{pre-l1})-(\ref{pre-l3}) in $\g$, the terms of the form $T ( \cdots T ( \cdots ) )$ cancel after using (\ref{lw}), and the remaining terms vanish by the representation axioms of $\h$ and the identities (\ref{comp-1})-(\ref{comp-9}). Alternatively, all nine identities follow from the fact that, for any subalgebra $K$ of a pre-Leibniz algebra $L$, the quotient $L/K$ is a representation of $K$; here $L = \g \ltimes_\lambda \h$, $K = Gr(T) \cong \h_T$ (Proposition \ref{prop-graph}), and $x \mapsto (x, 0) + K$ is an isomorphism $\g \cong L/K$ which carries the actions of $K$ to the above actions, since $(T(u), u) \ast_\lambda (x, 0) = \big( T(u) \ast x - T(u \ast^R x), 0 \big) + \big( T(u \ast^R x), u \ast^R x \big)$ and similarly on the right.
This proves that $\g$ is a representation of $\h_T$.
\end{proof}

Therefore, by the previous proposition, we can consider the cochain complex $\{ C^\ast_\pL (\h_T, \g), \delta_\pL \}$ of the pre-Leibniz algebra $\h_T$ with coefficients in the representation $\g$ of Proposition \ref{prop-rep-T}. In \cite{das-wrb}, the $0$-cochains are the elements $x_0 \in \mathfrak{g}$, which act through the pair of derivations $\big( [x_0, -], \rho^L (x_0, -) \big)$. Since left multiplications are not derivations of a pre-Leibniz algebra in general, we replace such pairs by the space $\mathrm{Der}(\g, \h)$ of pairs $(\phi, \psi)$, where $\phi$ is a derivation of $\g$ and $\psi$ is a derivation of $\h$ such that
\[
\psi (x \ast^L u) = \phi (x) \ast^L u + x \ast^L \psi (u), \qquad \psi (u \ast^R x) = \psi (u) \ast^R x + u \ast^R \phi (x), \qquad \forall x \in \g, ~ u \in \h, ~ \ast \in \{ \triangleleft, \triangleright \}.
\]
More precisely, $C^0_\pL (\h_T, \g) = \mathrm{Der}(\g, \h)$ and $C^n_\pL (\h_T, \g) = \mathrm{Hom}(\kk[C_n] \otimes \h^{\otimes n}, \g)$, for $n \geq 1$, and the differential $\delta_\pL : C^n_\pL (\h_T, \g) \rightarrow C^{n+1}_\pL (\h_T, \g)$ is given by
\[
\delta_\pL (\phi, \psi) = T \circ \psi - \phi \circ T, \qquad \forall (\phi, \psi) \in \mathrm{Der}(\g, \h),
\]
for $n = 0$, and by
\begin{equation}\label{delta-T}
\begin{split}
& (\delta_\pL f)([r]; u_1, \ldots, u_{n+1}) \\
&= \sum_{i=1}^n (-1)^{i+1} \Big( T(u_i) \ast_{R_i[r]} f (S_i[r]; u_1, \ldots, \widehat{u_i}, \ldots, u_{n+1}) - T \big( u_i \ast^R_{R_i[r]} f (S_i[r]; u_1, \ldots, \widehat{u_i}, \ldots, u_{n+1}) \big) \Big) \\
&\quad + (-1)^{n+1} \Big( f (S'[r]; u_1, \ldots, u_n) \ast_{R'[r]} T(u_{n+1}) - T \big( f (S'[r]; u_1, \ldots, u_n) \ast^L_{R'[r]} u_{n+1} \big) \Big) \\
&\quad + \sum_{1 \leq i < j \leq n+1} (-1)^i f \big( R_{ij}[r]; u_1, \ldots, \widehat{u_i}, \ldots, u_{j-1}, \\
&\hspace{7em} T(u_i) \ast^L_{S_{ij}[r]} u_j + u_i \ast^R_{S_{ij}[r]} T(u_j) + \lambda \, u_i \ast_{\h, S_{ij}[r]} u_j, u_{j+1}, \ldots, u_{n+1} \big),
\end{split}
\end{equation}
where the pairs $(R_i, S_i)$, $(R', S')$ and $(R_{ij}, S_{ij})$ are as in (\ref{pl-diff}), and, for $a \in \kk[C_2]$, the symbol $\ast_a$ stands for $\triangleleft$, $\triangleright$ or $\triangleleft + \triangleright$ according as $a = [1]$, $[2]$ or $[1] + [2]$ (and similarly for $\ast^L_a$, $\ast^R_a$ and $\ast_{\h, a}$), for $n \geq 1$. Note that $\delta_\pL \circ \delta_\pL = 0$ also on $C^0_\pL (\h_T, \g)$: for $(\phi, \psi) \in \mathrm{Der}(\g, \h)$, the map $(\mathrm{id}_\g - t \phi) \circ T \circ (\mathrm{id}_\h + t \psi) = T + t \, \delta_\pL (\phi, \psi)$ is a $\lambda$-weighted $\mathcal{O}$-operator modulo $t^2$, which means that $\delta_\pL (\phi, \psi)$ is a $1$-cocycle.

Let $Z^n_\pL (\h_T, \g)$ and $B^n_\pL (\h_T, \g)$ be the spaces of $n$-cocycles and $n$-coboundaries, respectively. Then we have $B^n_\pL (\h_T, \g) \subset Z^n_\pL (\h_T, \g)$. The corresponding quotients
\[
H^n_T (\h, \g) := \frac{Z^n_\pL (\h_T, \g)}{B^n_\pL (\h_T, \g)}, \qquad \text{for } n \geq 0,
\]
are called the cohomology groups of the $\lambda$-weighted $\mathcal{O}$-operator $T$.

\begin{remark}\label{rmk-dt}
It follows from the expressions (\ref{dT-explicit}) and (\ref{delta-T}) that $d_T f = (-1)^n \delta_\pL f$, $\forall f \in \mathrm{Hom}(\kk[C_n] \otimes \h^{\otimes n}, \g)$ with $n \geq 1$. Hence the cohomology groups induced by the differential $d_T$ are isomorphic to the cohomology groups of $T$.
\end{remark}

We now show that the cohomology of $T$ splits the cohomology of the corresponding weighted relative Rota-Baxter operator on Leibniz algebras \cite{das-wrb}.

\begin{prop}\label{prop-split}
Let $T : \h \rightarrow \g$ be a $\lambda$-weighted $\mathcal{O}$-operator. Then $\h_\Tot$ is a Leibniz $\g_\Tot$-representation and $T : \h_\Tot \rightarrow \g_\Tot$ is a $\lambda$-weighted relative Rota-Baxter operator in the sense of \cite{das-wrb}. Moreover, the maps $\Phi_n$ of Theorem \ref{thm-phi} define a morphism of cochain complexes from $\{ C^\ast_\pL (\h_T, \g), \delta_\pL \}$ to the cochain complex of the operator $T : \h_\Tot \rightarrow \g_\Tot$ given in \cite{das-wrb}. Hence there are morphisms $H^n_T (\h, \g) \rightarrow H^n_\Leib \big( (\h_\Tot)_T, \g_\Tot \big)$, for $n \geq 2$.
\end{prop}

\begin{proof}
The total Leibniz algebra of the pre-Leibniz algebra $\g \oplus \h$ with the operations (\ref{wsd-1}) has the bracket
\[
[(x, u), (y, v)] = \big( [x, y]_{\triangleleft, \triangleright}, ~ \rho^L (x, v) + \rho^R (u, y) + [u, v]_{\triangleleft_\h, \triangleright_\h} \big),
\]
where $\rho^L = \triangleleft^L + \triangleright^L$ and $\rho^R = \triangleleft^R + \triangleright^R$. Hence $\h_\Tot$ is a Leibniz $\g_\Tot$-representation. Adding the identities (\ref{lw}) for $\triangleleft$ and $\triangleright$, we get
\[
[T(u), T(v)]_{\triangleleft, \triangleright} = T \big( \rho^L (T(u), v) + \rho^R (u, T(v)) + \lambda [u, v]_{\triangleleft_\h, \triangleright_\h} \big), \qquad \forall u, v \in \h,
\]
so $T : \h_\Tot \rightarrow \g_\Tot$ is a $\lambda$-weighted relative Rota-Baxter operator. By (\ref{induced-pl}), we have $u \triangleleft_T v + u \triangleright_T v = [u, v]_T$, that is, $(\h_T)_\Tot = (\h_\Tot)_T$. Similarly, adding the actions of Proposition \ref{prop-rep-T}, we get
\[
u \triangleleft^L_T x + u \triangleright^L_T x = [T(u), x] - T \rho^R (u, x), \qquad x \triangleleft^R_T u + x \triangleright^R_T u = [x, T(u)] - T \rho^L (x, u),
\]
which are the actions of $(\h_\Tot)_T$ on $\g_\Tot$ given in \cite[Proposition 4.2]{das-wrb}. Therefore $\g_\Tot$, as a representation of $(\h_\Tot)_T$, is the total representation of the representation $\g$ of $\h_T$. The result now follows from Theorem \ref{thm-phi} applied to the pre-Leibniz algebra $\h_T$ and its representation $\g$.
\end{proof}

\section{Deformations of weighted \texorpdfstring{$\mathcal{O}$}{O}-operators}\label{sec-5}

In this section, following the classical approach of Gerstenhaber \cite{gers}, we study formal and finite order deformations of a $\lambda$-weighted $\mathcal{O}$-operator in terms of the cohomology introduced in the previous section. In particular, we find a sufficient condition for the rigidity of a $\lambda$-weighted $\mathcal{O}$-operator, and an obstruction class for the extensibility of a finite order deformation.

\medskip

\noindent \textbf{Formal deformations.} Let $\g, \h$ be two pre-Leibniz algebras and $\h$ be a pre-Leibniz $\g$-representation. We consider the spaces $\g[[t]]$ and $\h[[t]]$ of formal power series in $t$ with coefficients from $\g$ and $\h$, respectively. Note that $\g[[t]]$ and $\h[[t]]$ are both $\kk[[t]]$-modules. It is easy to see that the pre-Leibniz algebra structures on $\g$ and $\h$ can be extended by $\kk[[t]]$-bilinearity to pre-Leibniz algebra structures on $\g[[t]]$ and $\h[[t]]$, respectively. Moreover, the pre-Leibniz $\g$-representation on $\h$ also induces a pre-Leibniz $\g[[t]]$-representation on $\h[[t]]$.

\begin{defn}
Let $T : \h \rightarrow \g$ be a $\lambda$-weighted $\mathcal{O}$-operator. A formal deformation of $T$ consists of a formal sum $T_t = \sum_{i = 0}^\infty t^i T_i \in \mathrm{Hom}(\h, \g)[[t]]$ with $T_0 = T$, such that the $\kk[[t]]$-linear map $T_t : \h[[t]] \rightarrow \g[[t]]$ is a $\lambda$-weighted $\mathcal{O}$-operator (on $\h[[t]]$ over the pre-Leibniz algebra $\g[[t]]$).
\end{defn}

It follows from the above definition that $T_t$ is a formal deformation of $T$ if and only if
\begin{equation}\label{def-eqns}
\sum_{i+j = n} T_i (u) \ast T_j (v) = \sum_{i+j = n} T_i \big( T_j (u) \ast^L v + u \ast^R T_j (v) \big) + \lambda \, T_n (u \ast_\h v), \qquad \forall u, v \in \h, ~ n \geq 0, ~ \ast \in \{ \triangleleft, \triangleright \}.
\end{equation}
The above system of equations is called the deformation equations. Note that (\ref{def-eqns}) holds for $n = 0$ as $T_0 = T$ is a $\lambda$-weighted $\mathcal{O}$-operator. However, for $n = 1$, we get
\[
T(u) \ast T_1 (v) + T_1 (u) \ast T(v) = T \big( T_1 (u) \ast^L v + u \ast^R T_1 (v) \big) + T_1 \big( T(u) \ast^L v + u \ast^R T(v) + \lambda \, u \ast_\h v \big).
\]
This implies that $\delta_\pL (T_1) = 0$ (follows from (\ref{delta-T})). In other words, $T_1$ is a $1$-cocycle in the cohomology complex of $T$. This is called the infinitesimal of the deformation.

\medskip

\begin{defn}
Two deformations $T_t$ and $T'_t$ of a $\lambda$-weighted $\mathcal{O}$-operator $T$ are said to be equivalent if there exist a pair $(\phi, \psi) \in \mathrm{Der}(\g, \h)$ and linear maps $\phi_i \in \mathrm{Hom}(\g, \g)$, $\psi_i \in \mathrm{Hom}(\h, \h)$ for $i \geq 2$, such that the pair of maps
\[
\Big( \Phi_t = \mathrm{id}_\g + t \phi + \sum_{i \geq 2} t^i \phi_i : \g[[t]] \rightarrow \g[[t]], ~ \Psi_t = \mathrm{id}_\h + t \psi + \sum_{i \geq 2} t^i \psi_i : \h[[t]] \rightarrow \h[[t]] \Big)
\]
is a morphism of $\lambda$-weighted $\mathcal{O}$-operators from $T_t$ to $T'_t$.
\end{defn}

Therefore, if $T_t$ and $T'_t$ are equivalent, then we must have
\[
(\Phi_t \circ T_t)(u) = (T'_t \circ \Psi_t)(u), \qquad \Psi_t (x \ast^L u) = \Phi_t (x) \ast^L \Psi_t (u) \quad \text{and} \quad \Psi_t (u \ast^R x) = \Psi_t (u) \ast^R \Phi_t (x),
\]
$\forall x \in \g$, $u \in \h$ and $\ast \in \{ \triangleleft, \triangleright \}$. By expanding the identity $(\Phi_t \circ T_t)(u) = (T'_t \circ \Psi_t)(u)$ and equating coefficients of $t$ from both sides, we get
\[
T_1 (u) - T'_1 (u) = T (\psi (u)) - \phi (T(u)) = \delta_\pL (\phi, \psi)(u).
\]
As a summary of the above discussions, we get the following.

\begin{thm}\label{inf-1co}
Let $T : \h \rightarrow \g$ be a $\lambda$-weighted $\mathcal{O}$-operator. If $T_t = \sum_{i = 0}^\infty t^i T_i$ is a formal deformation of $T$, then $T_1$ is a $1$-cocycle in the cohomology complex of $T$. Moreover, the corresponding cohomology class depends only on the equivalence class of the deformation $T_t$.
\end{thm}

\begin{defn}
A $\lambda$-weighted $\mathcal{O}$-operator $T$ is said to be rigid if any formal deformation $T_t$ is equivalent to the undeformed one $T'_t = T$.
\end{defn}

Next, we introduce certain specific elements associated to $T$ whose importance will be clear in the next theorem.

\begin{defn}
A pair $(\phi, \psi) \in \mathrm{Der}(\g, \h)$ is said to be a Nijenhuis pair if
\[
\phi (x) \ast \phi (y) = 0, \qquad \psi (u) \ast_\h \psi (v) = 0, \qquad \phi (x) \ast^L \psi (u) = 0, \qquad \psi (u) \ast^R \phi (x) = 0,
\]
$\forall x, y \in \g$, $u, v \in \h$ and $\ast \in \{ \triangleleft, \triangleright \}$. We denote the set of all Nijenhuis pairs by $\mathrm{Nij}(\g, \h)$.
\end{defn}

We now have the following result.

\begin{thm}\label{thm-rigid}
Let $T : \h \rightarrow \g$ be a $\lambda$-weighted $\mathcal{O}$-operator. If $Z^1_\pL (\h_T, \g) = \delta_\pL (\mathrm{Nij}(\g, \h))$, then $T$ is rigid.
\end{thm}

\begin{proof}
Let $T_t = \sum_{i = 0}^\infty t^i T_i$ be a formal deformation of $T$. Then by Theorem \ref{inf-1co}, we get that $T_1 \in Z^1_\pL (\h_T, \g)$ is a $1$-cocycle in the cohomology complex of $T$. Hence by hypothesis, there exists a pair $(\phi, \psi) \in \mathrm{Nij}(\g, \h)$ such that $T_1 = \delta_\pL (\phi, \psi)$ (i.e. $T_1 (u) = T(\psi (u)) - \phi (T(u))$, $\forall u \in \h$). Set $\Phi_t = \mathrm{id}_\g + t \phi$ and $\Psi_t = \mathrm{id}_\h + t \psi$, and define $T'_t = \Phi_t \circ T_t \circ \Psi_t^{-1}$. Then $T'_t$ is a formal deformation of $T$. Since $(\phi, \psi) \in \mathrm{Nij}(\g, \h)$, it follows that $\Phi_t$ and $\Psi_t$ are pre-Leibniz algebra morphisms, and also $\Psi_t (x \ast^L u) = \Phi_t (x) \ast^L \Psi_t (u)$ and $\Psi_t (u \ast^R x) = \Psi_t (u) \ast^R \Phi_t (x)$ hold. Hence $T'_t$ is equivalent to $T_t$. We also have
\begin{align*}
T'_t (u) &= (\mathrm{id}_\g + t \phi) \circ \Big( \sum_{i = 0}^\infty t^i T_i \Big) \big( u - t \psi (u) + \text{powers of } t^{\geq 2} \big) \\
&= (\mathrm{id}_\g + t \phi) \big( T(u) - t \, T(\psi (u)) + t \, T_1 (u) + \text{powers of } t^{\geq 2} \big) \\
&= T(u) + t \big( T_1 (u) - T(\psi (u)) + \phi (T(u)) \big) + \text{powers of } t^{\geq 2}.
\end{align*}
It follows from the above expression that the coefficient of $t$ vanishes. By repeating the above argument, we get that $T_t$ is equivalent to $T$.
\end{proof}

\noindent \textbf{Finite order deformations.} Let $T : \h \rightarrow \g$ be a $\lambda$-weighted $\mathcal{O}$-operator. For a fixed $N \in \mathbb{N}$, consider the spaces $\g[[t]]/(t^{N+1})$ and $\h[[t]]/(t^{N+1})$. Both of them are $\kk[[t]]/(t^{N+1})$-modules. Moreover, by $\kk[[t]]/(t^{N+1})$-bilinearity, the pre-Leibniz structures on $\g$ and $\h$ can be extended to $\g[[t]]/(t^{N+1})$ and $\h[[t]]/(t^{N+1})$, respectively. Further, $\h[[t]]/(t^{N+1})$ is a pre-Leibniz $\g[[t]]/(t^{N+1})$-representation.

\begin{defn}
An order $N$ deformation of $T$ is given by a finite polynomial of the form $T_t^N = \sum_{i = 0}^N t^i T_i \in \mathrm{Hom}(\h, \g)[[t]]/(t^{N+1})$ with $T_0 = T$, such that the $\kk[[t]]/(t^{N+1})$-linear map $T^N_t : \h[[t]]/(t^{N+1}) \rightarrow \g[[t]]/(t^{N+1})$ is a $\lambda$-weighted $\mathcal{O}$-operator.
\end{defn}

Thus, $T_t^N = \sum_{i = 0}^N t^i T_i$ is an order $N$ deformation of $T$ if and only if the identity (\ref{def-eqns}) holds for $n = 0, 1, \ldots, N$. These can be equivalently written as
\begin{equation}\label{fin-eq}
d_T (T_n) = - \frac{1}{2} \sum_{ \substack{i+j = n\\ i, j \geq 1}} \llbracket T_i, T_j \rrbracket, \qquad \text{for } n = 1, \ldots, N.
\end{equation}

\begin{defn}
(Obstruction cochain) Let $T_t^N = \sum_{i = 0}^N t^i T_i$ be an order $N$ deformation of $T$. We define a $2$-cochain $Ob_{T_t^N} \in C^2_\pL (\h_T, \g)$, called the obstruction cochain, by
\[
Ob_{T_t^N} = - \frac{1}{2} \sum_{\substack{i+j = N+1\\ i, j \geq 1}} \llbracket T_i, T_j \rrbracket.
\]
\end{defn}

\begin{prop}
The obstruction cochain $Ob_{T_t^N}$ is a $2$-cocycle in the cohomology complex of $T$, i.e. $\delta_\pL (Ob_{T_t^N}) = 0$ (equivalently, $d_T (Ob_{T_t^N}) = 0$).
\end{prop}

\begin{proof}
We have
\begin{align*}
d_T \Big( - \frac{1}{2} \sum_{ \substack{ i+j = N+1\\ i, j \geq 1}} \llbracket T_i, T_j \rrbracket \Big)
&= - \frac{1}{2} \sum_{\substack{  i+j = N+1\\ i, j \geq 1}} \big( d \llbracket T_i, T_j \rrbracket + \llbracket T, \llbracket T_i, T_j \rrbracket \rrbracket \big) \\
&= - \frac{1}{2} \sum_{\substack{  i+j = N+1\\ i, j \geq 1}} \big( \llbracket dT_i, T_j \rrbracket - \llbracket T_i, dT_j \rrbracket + \llbracket \llbracket T, T_i \rrbracket, T_j \rrbracket - \llbracket T_i, \llbracket T, T_j \rrbracket \rrbracket \big) \\
&= - \frac{1}{2} \sum_{\substack{  i+j = N+1\\ i, j \geq 1}} \big( \llbracket d_T T_i, T_j \rrbracket -  \llbracket T_i, d_T T_j \rrbracket \big)  \\
&= \frac{1}{4} \sum_{\substack{  i_1 +  i_2 + j = N+1\\ i_1, i_2, j \geq 1}} \llbracket  \llbracket T_{i_1}, T_{i_2} \rrbracket, T_j \rrbracket - \frac{1}{4} \sum_{\substack{  i + j_1 + j_2 = N+1\\ i, j_1, j_2 \geq 1}} \llbracket T_i, \llbracket T_{j_1}, T_{j_2} \rrbracket \rrbracket \qquad (\text{from } (\ref{fin-eq})) \\
&= \frac{1}{2} \sum_{\substack{  i + j + k = N+1\\ i, j, k \geq 1}} \llbracket  \llbracket T_{i}, T_{j} \rrbracket, T_k \rrbracket = 0.
\end{align*}
This completes the proof.
\end{proof}

It follows from the above result that we obtain a cohomology class $[Ob_{T_t^N}] \in H^2_T (\h, \g)$ associated to the order $N$ deformation $T^N_t$. This cohomology class is called the obstruction class.

\begin{defn}
An order $N$ deformation $T_t^N = \sum_{i = 0}^N t^i T_i$ is said to be extensible if there exists a linear map $T_{N+1} : \h \rightarrow \g$ which makes $T^{N+1}_t = T^N_t + t^{N+1} T_{N+1}$ into an order $N+1$ deformation.
\end{defn}

The following result gives a necessary and sufficient condition for the extensibility of a finite order deformation in terms of its obstruction class.

\begin{thm}
An order $N$ deformation $T_t^N$ is extensible if and only if the corresponding obstruction class $[Ob_{T_t^N}]$ vanishes.
\end{thm}

\begin{proof}
Let $T_t^N = \sum_{i = 0}^N t^i T_i$ be an extensible order $N$ deformation. Since there exists a linear map $T_{N+1}$ which makes $T^{N+1}_t = \sum_{i = 0}^{N+1} t^i T_i$ into an order $N+1$ deformation, we have that $Ob_{T_t^N} = d_T (T_{N+1})$. In other words, $Ob_{T_t^N}$ is a coboundary. Hence the cohomology class $[Ob_{T_t^N}]$ vanishes. The converse part follows by a similar argument.
\end{proof}

\section{Post-pre-Leibniz algebras}\label{sec-6}

In this section, we introduce post-pre-Leibniz algebras as the algebraic structures underlying $\lambda$-weighted $\mathcal{O}$-operators on pre-Leibniz algebras. They are to pre-Leibniz algebras what post-Leibniz algebras \cite{das-wrb} are to Leibniz algebras. In the terminology of \cite{bai-operad}, where the splittings of an operad by Rota-Baxter operators of weight zero and of nonzero weight are called its bisuccessor and its trisuccessor, post-pre-Leibniz algebras are the trisuccessors of pre-Leibniz algebras; this explains the name.

\begin{defn}\label{defn-ppl}
A post-pre-Leibniz algebra is a vector space $\h$ together with six bilinear operations $\triangleleft_1, \triangleleft_2, \triangleleft_3, \triangleright_1, \triangleright_2, \triangleright_3 : \h \otimes \h \rightarrow \h$ such that, setting
\[
x \Lt{\star} y = x \Lt{1} y + x \Lt{2} y + x \Lt{3} y, \qquad x \Rt{\star} y = x \Rt{1} y + x \Rt{2} y + x \Rt{3} y,
\]
the following twenty-one identities hold $\forall x, y, z \in \h$. They are all obtained from (\ref{pre-l1})-(\ref{pre-l3}) by a single rule, described in Remark \ref{rmk-color}, and are grouped accordingly:
\begin{align}
\shortintertext{\textit{Type} $S = \{ x \}$:}
x \Lt{1} (y \Lt{\star} z + y \Rt{\star} z) &= (x \Lt{1} y) \Lt{1} z + y \Rt{2} (x \Lt{1} z) \label{ppl-x-1}\\
x \Rt{1} (y \Lt{\star} z) &= (x \Rt{1} y) \Lt{1} z + y \Lt{2} (x \Lt{1} z + x \Rt{1} z) \label{ppl-x-2}\\
x \Rt{1} (y \Rt{\star} z) &= (x \Lt{1} y + x \Rt{1} y) \Rt{1} z + y \Rt{2} (x \Rt{1} z) \label{ppl-x-3}\\
\shortintertext{\textit{Type} $S = \{ y \}$:}
x \Lt{2} (y \Lt{1} z + y \Rt{1} z) &= (x \Lt{2} y) \Lt{1} z + y \Rt{1} (x \Lt{\star} z) \label{ppl-y-1}\\
x \Rt{2} (y \Lt{1} z) &= (x \Rt{2} y) \Lt{1} z + y \Lt{1} (x \Lt{\star} z + x \Rt{\star} z) \label{ppl-y-2}\\
x \Rt{2} (y \Rt{1} z) &= (x \Lt{2} y + x \Rt{2} y) \Rt{1} z + y \Rt{1} (x \Rt{\star} z) \label{ppl-y-3}\\
\shortintertext{\textit{Type} $S = \{ z \}$:}
x \Lt{2} (y \Lt{2} z + y \Rt{2} z) &= (x \Lt{\star} y) \Lt{2} z + y \Rt{2} (x \Lt{2} z) \label{ppl-z-1}\\
x \Rt{2} (y \Lt{2} z) &= (x \Rt{\star} y) \Lt{2} z + y \Lt{2} (x \Lt{2} z + x \Rt{2} z) \label{ppl-z-2}\\
x \Rt{2} (y \Rt{2} z) &= (x \Lt{\star} y + x \Rt{\star} y) \Rt{2} z + y \Rt{2} (x \Rt{2} z) \label{ppl-z-3}\\
\shortintertext{\textit{Type} $S = \{ y, z \}$:}
x \Lt{2} (y \Lt{3} z + y \Rt{3} z) &= (x \Lt{2} y) \Lt{3} z + y \Rt{3} (x \Lt{2} z) \label{ppl-yz-1}\\
x \Rt{2} (y \Lt{3} z) &= (x \Rt{2} y) \Lt{3} z + y \Lt{3} (x \Lt{2} z + x \Rt{2} z) \label{ppl-yz-2}\\
x \Rt{2} (y \Rt{3} z) &= (x \Lt{2} y + x \Rt{2} y) \Rt{3} z + y \Rt{3} (x \Rt{2} z) \label{ppl-yz-3}\\
\shortintertext{\textit{Type} $S = \{ x, z \}$:}
x \Lt{3} (y \Lt{2} z + y \Rt{2} z) &= (x \Lt{1} y) \Lt{3} z + y \Rt{2} (x \Lt{3} z) \label{ppl-xz-1}\\
x \Rt{3} (y \Lt{2} z) &= (x \Rt{1} y) \Lt{3} z + y \Lt{2} (x \Lt{3} z + x \Rt{3} z) \label{ppl-xz-2}\\
x \Rt{3} (y \Rt{2} z) &= (x \Lt{1} y + x \Rt{1} y) \Rt{3} z + y \Rt{2} (x \Rt{3} z) \label{ppl-xz-3}\\
\shortintertext{\textit{Type} $S = \{ x, y \}$:}
x \Lt{3} (y \Lt{1} z + y \Rt{1} z) &= (x \Lt{3} y) \Lt{1} z + y \Rt{3} (x \Lt{1} z) \label{ppl-xy-1}\\
x \Rt{3} (y \Lt{1} z) &= (x \Rt{3} y) \Lt{1} z + y \Lt{3} (x \Lt{1} z + x \Rt{1} z) \label{ppl-xy-2}\\
x \Rt{3} (y \Rt{1} z) &= (x \Lt{3} y + x \Rt{3} y) \Rt{1} z + y \Rt{3} (x \Rt{1} z) \label{ppl-xy-3}\\
\shortintertext{\textit{Type} $S = \{ x, y, z \}$:}
x \Lt{3} (y \Lt{3} z + y \Rt{3} z) &= (x \Lt{3} y) \Lt{3} z + y \Rt{3} (x \Lt{3} z) \label{ppl-xyz-1}\\
x \Rt{3} (y \Lt{3} z) &= (x \Rt{3} y) \Lt{3} z + y \Lt{3} (x \Lt{3} z + x \Rt{3} z) \label{ppl-xyz-2}\\
x \Rt{3} (y \Rt{3} z) &= (x \Lt{3} y + x \Rt{3} y) \Rt{3} z + y \Rt{3} (x \Rt{3} z) \label{ppl-xyz-3}
\end{align}
\end{defn}

\begin{remark}\label{rmk-color}
The identities of Definition \ref{defn-ppl} can be described as follows. Consider (\ref{pre-l1})-(\ref{pre-l3}) and choose a nonempty subset $S \subset \{ x, y, z \}$. The variables in $S$ are said to be of type $\h$, and the other ones of type $\g$. Each product of two factors is replaced according to the types of these factors:
\[
\renewcommand{\arraystretch}{1.3}
\begin{array}{c|cccc}
\text{types of the factors} & (\g, \g) & (\g, \h) & (\h, \g) & (\h, \h) \\ \hline
\text{operation} & \ast_\star & \ast_2 & \ast_1 & \ast_3 \\
\text{type of the product} & \g & \h & \h & \h
\end{array}
\]
for $\ast \in \{ \triangleleft, \triangleright \}$. In other words, a product is of type $\g$ only when both factors are of type $\g$. For instance, for $S = \{ y \}$, the product $y \triangleleft z$ in (\ref{pre-l1}) becomes $y \Lt{1} z$, of type $\h$, so that $x \triangleleft (y \triangleleft z)$ becomes $x \Lt{2} (y \Lt{1} z)$, while $x \triangleleft z$ becomes $x \Lt{\star} z$; this gives (\ref{ppl-y-1}). The $7 \times 3 = 21$ resulting identities are exactly the ones of Definition \ref{defn-ppl}, listed by groups of three (coming from (\ref{pre-l1}), (\ref{pre-l2}), (\ref{pre-l3}) respectively) for $S = \{x\}, \{y\}, \{z\}, \{y,z\}, \{x,z\}, \{x,y\}, \{x,y,z\}$. The three identities with $S = \{ x, y, z \}$ say that $(\h, \triangleleft_3, \triangleright_3)$ is a pre-Leibniz algebra. The case $S = \emptyset$ is not needed: the corresponding identities follow from the other ones (Proposition \ref{prop-ppl-tot}). The types correspond to the two summands of the pre-Leibniz algebra $\h_\Tot \oplus \h$ in Proposition \ref{prop-converse}.
\end{remark}

As for post-Leibniz algebras, the sum of the operations of a post-pre-Leibniz algebra gives a pre-Leibniz algebra.

\begin{prop}\label{prop-ppl-tot}
Let $(\h, \triangleleft_1, \triangleleft_2, \triangleleft_3, \triangleright_1, \triangleright_2, \triangleright_3)$ be a post-pre-Leibniz algebra. Then $(\h, \triangleleft_\star, \triangleright_\star)$ is a pre-Leibniz algebra, called the total pre-Leibniz algebra and denoted by $\h_\Tot$.
\end{prop}

\begin{proof}
By adding the left hand sides of the identities (\ref{ppl-x-1}), (\ref{ppl-y-1}), (\ref{ppl-z-1}), (\ref{ppl-yz-1}), (\ref{ppl-xz-1}), (\ref{ppl-xy-1}) and (\ref{ppl-xyz-1}), we simply get $x \Lt{\star} (y \Lt{\star} z + y \Rt{\star} z)$. On the other hand, by adding the right hand sides of these identities, we get $(x \Lt{\star} y) \Lt{\star} z + y \Rt{\star} (x \Lt{\star} z)$. Therefore, the identity (\ref{pre-l1}) holds for the operations $\triangleleft_\star$ and $\triangleright_\star$. Similarly, by adding the second (resp. third) identities of each type, we get the identity (\ref{pre-l2}) (resp. (\ref{pre-l3})).
\end{proof}

The next result shows that post-pre-Leibniz algebras are the structures induced by $\lambda$-weighted $\mathcal{O}$-operators.

\begin{prop}\label{prop-T-ppl}
Let $\h$ be a pre-Leibniz $\g$-representation and $T : \h \rightarrow \g$ be a $\lambda$-weighted $\mathcal{O}$-operator. Then $\h$ carries a post-pre-Leibniz algebra structure given by
\[
u \ast_1 v = u \ast^R T(v), \qquad u \ast_2 v = T(u) \ast^L v, \qquad u \ast_3 v = \lambda \, u \ast_\h v, \qquad \forall u, v \in \h, ~ \ast \in \{ \triangleleft, \triangleright \}.
\]
Its total pre-Leibniz algebra is the induced pre-Leibniz algebra $\h_T$.
\end{prop}

\begin{proof}
Let $u, v, w \in \h$. Recall that $u \ast_\star v = u \ast_T v$ and $T(u \ast_T v) = T(u) \ast T(v)$. We verify the first identity of each type in Definition \ref{defn-ppl}.

\noindent \emph{Type $S = \{x\}$.} We have
\begin{align*}
u \Lt{1} \big( v \Lt{\star} w + v \Rt{\star} w \big) &= u \triangleleft^R T \big( v \triangleleft_T w + v \triangleright_T w \big) \\
&= u \triangleleft^R \big( T(v) \triangleleft T(w) + T(v) \triangleright T(w) \big), \\
(u \Lt{1} v) \Lt{1} w + v \Rt{2} (u \Lt{1} w) &= \big( u \triangleleft^R T(v) \big) \triangleleft^R T(w) + T(v) \triangleright^L \big( u \triangleleft^R T(w) \big).
\end{align*}
The right hand sides are equal by the identity (\ref{pre-l1}) for the representation $\h$ of $\g$, with arguments $(u, T(v), T(w))$. Hence the identity holds.

\noindent \emph{Type $S = \{y\}$.} We have
\begin{align*}
u \Lt{2} \big( v \Lt{1} w + v \Rt{1} w \big) &= T(u) \triangleleft^L \big( v \triangleleft^R T(w) + v \triangleright^R T(w) \big), \\
(u \Lt{2} v) \Lt{1} w + v \Rt{1} (u \Lt{\star} w) &= \big( T(u) \triangleleft^L v \big) \triangleleft^R T(w) + v \triangleright^R \big( T(u) \triangleleft T(w) \big),
\end{align*}
which are equal by (\ref{pre-l1}) for the representation $\h$, with arguments $(T(u), v, T(w))$.

\noindent \emph{Type $S = \{z\}$.} We have
\begin{align*}
u \Lt{2} \big( v \Lt{2} w + v \Rt{2} w \big) &= T(u) \triangleleft^L \big( T(v) \triangleleft^L w + T(v) \triangleright^L w \big), \\
(u \Lt{\star} v) \Lt{2} w + v \Rt{2} (u \Lt{2} w) &= \big( T(u) \triangleleft T(v) \big) \triangleleft^L w + T(v) \triangleright^L \big( T(u) \triangleleft^L w \big),
\end{align*}
which are equal by (\ref{pre-l1}) for the representation $\h$, with arguments $(T(u), T(v), w)$.

\noindent \emph{Type $S = \{y, z\}$.} We have
\begin{align*}
u \Lt{2} \big( v \Lt{3} w + v \Rt{3} w \big) &= \lambda \, T(u) \triangleleft^L \big( v \triangleleft_\h w + v \triangleright_\h w \big), \\
(u \Lt{2} v) \Lt{3} w + v \Rt{3} (u \Lt{2} w) &= \lambda \big( T(u) \triangleleft^L v \big) \triangleleft_\h w + \lambda \, v \triangleright_\h \big( T(u) \triangleleft^L w \big),
\end{align*}
which are equal by (\ref{comp-1}) with $x = T(u) \in \g$.

\noindent \emph{Type $S = \{x, z\}$.} We have
\begin{align*}
u \Lt{3} \big( v \Lt{2} w + v \Rt{2} w \big) &= \lambda \, u \triangleleft_\h \big( T(v) \triangleleft^L w + T(v) \triangleright^L w \big), \\
(u \Lt{1} v) \Lt{3} w + v \Rt{2} (u \Lt{3} w) &= \lambda \big( u \triangleleft^R T(v) \big) \triangleleft_\h w + \lambda \, T(v) \triangleright^L \big( u \triangleleft_\h w \big),
\end{align*}
which are equal by (\ref{comp-4}) with $y = T(v) \in \g$.

\noindent \emph{Type $S = \{x, y\}$.} We have
\begin{align*}
u \Lt{3} \big( v \Lt{1} w + v \Rt{1} w \big) &= \lambda \, u \triangleleft_\h \big( v \triangleleft^R T(w) + v \triangleright^R T(w) \big), \\
(u \Lt{3} v) \Lt{1} w + v \Rt{3} (u \Lt{1} w) &= \lambda \big( u \triangleleft_\h v \big) \triangleleft^R T(w) + \lambda \, v \triangleright_\h \big( u \triangleleft^R T(w) \big),
\end{align*}
which are equal by (\ref{comp-7}) with $z = T(w) \in \g$.

\noindent \emph{Type $S = \{x, y, z\}$.} The identity reads $\lambda^2 \big( u \triangleleft_\h ( v \triangleleft_\h w + v \triangleright_\h w ) - ( u \triangleleft_\h v ) \triangleleft_\h w - v \triangleright_\h ( u \triangleleft_\h w ) \big) = 0$, which is (\ref{pre-l1}) in $\h$.

The second and the third identities of each type are obtained in the same way, using (\ref{pre-l2}) and (\ref{pre-l3}) instead of (\ref{pre-l1}). Finally, $u \ast_\star v = u \ast_T v$ by (\ref{induced-pl}). This completes the proof.
\end{proof}

Conversely, any post-pre-Leibniz algebra is induced by a $\lambda$-weighted $\mathcal{O}$-operator.

\begin{prop}\label{prop-converse}
Let $(\h, \triangleleft_1, \triangleleft_2, \triangleleft_3, \triangleright_1, \triangleright_2, \triangleright_3)$ be a post-pre-Leibniz algebra and $\lambda \neq 0$. Then the pre-Leibniz algebra $(\h, \tfrac{1}{\lambda} \triangleleft_3, \tfrac{1}{\lambda} \triangleright_3)$ is a pre-Leibniz $\h_\Tot$-representation with the actions
\[
x \ast^L u = x \ast_2 u, \qquad u \ast^R x = u \ast_1 x, \qquad \forall x \in \h_\Tot, ~ u \in \h, ~ \ast \in \{ \triangleleft, \triangleright \}.
\]
Moreover, the identity map $\mathrm{id} : \h \rightarrow \h_\Tot$ is a $\lambda$-weighted $\mathcal{O}$-operator, and the induced post-pre-Leibniz algebra structure on $\h$ is the given one.
\end{prop}

\begin{proof}
Write $u \ast_\h v := \frac{1}{\lambda} \, u \ast_3 v$. For $\lambda = 1$, the operations (\ref{wsd-1}) on $\h_\Tot \oplus \h$ are
\[
(x, u) \ast (y, v) = \big( x \ast_\star y, ~ x \ast_2 v + u \ast_1 y + u \ast_3 v \big).
\]
Hence the rules of Remark \ref{rmk-color} are exactly the products of $\h_\Tot \oplus \h$, with type $\g$ corresponding to $\h_\Tot \oplus 0$ and type $\h$ to $0 \oplus \h$. The identities (\ref{pre-l1})-(\ref{pre-l3}) for $\h_\Tot \oplus \h$ with all arguments in $\h_\Tot$ hold by Proposition \ref{prop-ppl-tot}, and the other ones are the twenty-one identities of Definition \ref{defn-ppl}. Hence $\h_\Tot \oplus \h$ is a pre-Leibniz algebra, i.e. $(\h, \triangleleft_3, \triangleright_3)$ is a pre-Leibniz $\h_\Tot$-representation. For an arbitrary $\lambda \neq 0$, the identities (\ref{comp-1})-(\ref{comp-9}) are linear in the operations of $\h$, and the identities (\ref{pre-l1})-(\ref{pre-l3}) of $\h$ are quadratic in them. Hence they still hold for $\frac{1}{\lambda} \triangleleft_3$ and $\frac{1}{\lambda} \triangleright_3$, and $(\h, \frac{1}{\lambda} \triangleleft_3, \frac{1}{\lambda} \triangleright_3)$ is also a pre-Leibniz $\h_\Tot$-representation. Next,
\[
\mathrm{id}(u) \ast_\star \mathrm{id}(v) = \mathrm{id} \big( u \ast_2 v + u \ast_1 v + \lambda \, u \ast_\h v \big), \qquad \forall u, v \in \h,
\]
which is (\ref{lw}) for $T = \mathrm{id}$. Finally, by Proposition \ref{prop-T-ppl}, the induced operations are $u \ast^R \mathrm{id}(v) = u \ast_1 v$, $\mathrm{id}(u) \ast^L v = u \ast_2 v$ and $\lambda \, u \ast_\h v = u \ast_3 v$, which is the given structure.
\end{proof}

Given a pre-Leibniz algebra, the following result gives a necessary and sufficient condition for the existence of a compatible post-pre-Leibniz algebra structure.

\begin{prop}
Let $(\g, \triangleleft, \triangleright)$ be a pre-Leibniz algebra and $\lambda \neq 0$. Then there is a compatible post-pre-Leibniz algebra structure on $\g$ (i.e. one whose total pre-Leibniz algebra is $(\g, \triangleleft, \triangleright)$) if and only if there exist a pre-Leibniz $\g$-representation $\h$ and an invertible $\lambda$-weighted $\mathcal{O}$-operator $T : \h \rightarrow \g$.
\end{prop}

\begin{proof}
If $\g$ has a compatible post-pre-Leibniz algebra structure, then by Proposition \ref{prop-converse}, the identity map from $(\g, \frac{1}{\lambda} \triangleleft_3, \frac{1}{\lambda} \triangleright_3)$ to $\g_\Tot = (\g, \triangleleft, \triangleright)$ is an invertible $\lambda$-weighted $\mathcal{O}$-operator. Conversely, let $T : \h \rightarrow \g$ be an invertible $\lambda$-weighted $\mathcal{O}$-operator. By Proposition \ref{prop-T-ppl}, $\h$ carries a post-pre-Leibniz algebra structure, which we transport to $\g$ via $T$: $\forall x, y \in \g$,
\[
x \ast_1 y = T \big( T^{-1}(x) \ast^R y \big), \quad x \ast_2 y = T \big( x \ast^L T^{-1}(y) \big), \quad x \ast_3 y = \lambda \, T \big( T^{-1}(x) \ast_\h T^{-1}(y) \big).
\]
Then, by (\ref{lw}),
\begin{align*}
x \ast_1 y + x \ast_2 y + x \ast_3 y &= T \big( T^{-1}(x) \ast^R T T^{-1}(y) + T T^{-1}(x) \ast^L T^{-1}(y) + \lambda \, T^{-1}(x) \ast_\h T^{-1}(y) \big) \\
&= T T^{-1}(x) \ast T T^{-1}(y) = x \ast y.
\end{align*}
Hence the transported structure is compatible with $(\g, \triangleleft, \triangleright)$.
\end{proof}

\begin{remark}\label{rmk-ppl-pl}
Let $(\h, \triangleleft_i, \triangleright_i)_{i = 1,2,3}$ be a post-pre-Leibniz algebra. Then $(\h, \Lt{1} + \Rt{1}, \Lt{2} + \Rt{2}, \Lt{3} + \Rt{3})$ is a post-Leibniz algebra in the sense of \cite{das-wrb}. When $\triangleleft_3 = \triangleright_3 = 0$ (e.g. for $\lambda = 0$ in Proposition \ref{prop-T-ppl}), we recover the Leibniz-quadri-algebras of \cite{sun-guo}, which are analogues of the quadri-algebras of \cite{aguiar-loday}.
\end{remark}

\begin{exam}
Let $\g$ be a pre-Leibniz algebra and $\g = \g_1 \oplus \g_2$ be a decomposition into two pre-Leibniz subalgebras. The projection $P_1$ onto $\g_1$ is a $(-1)$-weighted Rota-Baxter operator, since $P_1 ( P_1(x) \ast y + x \ast P_1(y) - x \ast y ) = P_1 ( x_1 \ast y_1 - x_2 \ast y_2 ) = P_1(x) \ast P_1(y)$ for $x = x_1 + x_2$ and $y = y_1 + y_2$ with $x_i, y_i \in \g_i$; so $- P_1$ is a $1$-weighted Rota-Baxter operator by Example \ref{exam-wrb} (iv). By Proposition \ref{prop-T-ppl} (applied to the adjoint representation), $\g$ carries the post-pre-Leibniz algebra structure
\[
x \ast_1 y = - x \ast P_1 (y), \qquad x \ast_2 y = - P_1 (x) \ast y, \qquad x \ast_3 y = x \ast y, \qquad \forall x, y \in \g,
\]
whose total pre-Leibniz algebra is given by (with $x = x_1 + x_2$ and $y = y_1 + y_2$)
\[
x \ast_\star y = x \ast y - P_1 (x) \ast y - x \ast P_1 (y) = x_2 \ast y_2 - x_1 \ast y_1,
\]
i.e. the direct product $\g_1^{\mathrm{op}} \times \g_2$, where $\g_1^{\mathrm{op}}$ denotes $\g_1$ with the operations multiplied by $-1$. For instance, one can take $\g = A \times \mathfrak{s}$, where $A$ is the algebra of $2 \times 2$ upper triangular matrices viewed as a pre-Leibniz algebra as in Example \ref{exam-pre} (ii) and $\mathfrak{s} = \mathfrak{sl}_2$ viewed as in Example \ref{exam-pre} (i), with $\g_1 = \kk E_{12} \times \mathfrak{n}^+$ and $\g_2 = (\kk E_{11} \oplus \kk E_{22}) \times \mathfrak{b}^-$. Indeed, $\kk E_{12}$ and $\kk E_{11} \oplus \kk E_{22}$ are subalgebras of $A$ with $A = \kk E_{12} \oplus (\kk E_{11} \oplus \kk E_{22})$, and $\mathfrak{n}^+ = \kk e$ and $\mathfrak{b}^- = \kk h \oplus \kk f$ are subalgebras of $\mathfrak{sl}_2$ with $\mathfrak{sl}_2 = \mathfrak{n}^+ \oplus \mathfrak{b}^-$.
\end{exam}

%\begin{remark}
%In \cite[Section 6]{das-pre}, homotopy pre-Leibniz algebras (pre-Leibniz$_\infty$-algebras) are introduced and Rota-Baxter operators on Leibniz$_\infty$-algebras are shown to induce pre-Leibniz$_\infty$-algebras. It would be interesting to define weighted $\mathcal{O}$-operators on pre-Leibniz$_\infty$-algebras and the corresponding homotopy post-pre-Leibniz algebras, and to find a Maurer-Cartan characterization (and cohomology) of post-pre-Leibniz algebras splitting the one of pre-Leibniz algebras \cite[Theorem 3.7]{das-pre}.
%\end{remark}

\medskip

\section*{Declarations}

\textbf{Funding.} The author did not receive support from any organization for the submitted work.

\medskip

\textbf{Conflict of interest.} The author declares that he has no conflict of interest.


\begin{thebibliography}{99}

%\setstretch{1.1}
\bibitem{aguiar-loday} M. Aguiar and J.-L. Loday, Quadri-algebras, {\em J. Pure Appl. Algebra} {\bf 191} (2004), 205--221.

\bibitem{bai-operad} C. Bai, O. Bellier, L. Guo and X. Ni, Splitting of operations, Manin products, and Rota-Baxter operators, {\em Int. Math. Res. Not. IMRN} {\bf 2013} (2013), no. 3, 485--524.

\bibitem{bgn2} C. Bai, L. Guo and X. Ni, $\mathcal{O}$-operators on associative algebras and associative Yang-Baxter equations, {\em Pacific J. Math.} {\bf 256} (2012), no. 2, 257--289.

\bibitem{bgn} C. Bai, L. Guo and X. Ni, Nonabelian generalized Lax pairs, the classical Yang-Baxter equation and PostLie algebras, {\em Comm. Math. Phys.} {\bf 297} (2010), 553--596.

\bibitem{bgst} C. Bai, L. Guo, Y. Sheng and R. Tang, Post-groups, (Lie-)Butcher groups and the Yang-Baxter equation, {\em Math. Ann.} {\bf 388} (2024), no. 3, 3127--3167.

\bibitem{bala-operad} D. Balavoine, Deformations of algebras over a quadratic operad, in: {\em Operads: Proceedings of Renaissance Conferences (Hartford, CT/Luminy, 1995)}, Contemp. Math. {\bf 202}, Amer. Math. Soc., Providence, RI, 1997, 207--234.

\bibitem{bb-lya} K. Basdouri and S. Benabdelhafidh, Cohomology and deformations of modified Rota-Baxter Lie-Yamaguti algebras, {\em C. R. Math. Acad. Sci. Paris} {\bf 363} (2025), 641--662.

\bibitem{bbt-jj} I. Basdouri, S. Benabdelhafidh and W. Teng, Cohomology of modified $\lambda$-differential Jacobi-Jordan algebras and its applications, {\em Acta Comment. Univ. Tartu. Math.} {\bf 28} (2024), no. 2, 215--232.

\bibitem{baxter} G. Baxter, An analytic problem whose solution follows from a simple algebraic identity, {\em Pacific J. Math.} {\bf 10} (1960), 731--742.

\bibitem{b-mor}
S. Benabdelhafidh,
Morphisms cohomology and deformations of modified
$\lambda$-differential algebras,
{\em J. Algebra Appl.} {\bf 26} (2027), 2750006.
\href{https://doi.org/10.1142/S021949882750006X}
{doi:10.1142/S021949882750006X}.

\bibitem{bloh} A. Bloh, A generalization of the concept of a Lie algebra, {\em Dokl. Akad. Nauk SSSR} {\bf 165} (1965), no. 3, 471--473.

\bibitem{cartier} P. Cartier, On the structure of free Baxter algebras, {\em Adv. Math.} {\bf 9} (1972), 253--265.

\bibitem{connes} A. Connes and D. Kreimer, Renormalization in quantum field theory and the Riemann-Hilbert problem. I. The Hopf algebra structure of graphs and the main theorem, {\em Comm. Math. Phys.} {\bf 210} (2000), 249--273.

\bibitem{das-weighted} A. Das, Cohomology and deformations of weighted Rota-Baxter operators, {\em J. Math. Phys.} {\bf 63} (2022), 091703.

\bibitem{das-rota} A. Das, Deformations of associative Rota-Baxter operators, {\em J. Algebra} {\bf 560} (2020), 144--180.

\bibitem{das-pre} A. Das, Pre-Leibniz algebras, {\em Comm. Algebra} {\bf 52} (2024), no. 8, 3383--3399.

\bibitem{das-wrb} A. Das, Relative Rota-Baxter operators of arbitrary weight on Leibniz algebras and post-Leibniz algebra structures, {\em Publ. Math. Debrecen} {\bf 103} (2023), no. 3--4, 385--406.

\bibitem{ebrahimi} K. Ebrahimi-Fard, Loday-type algebras and the Rota-Baxter relation, {\em Lett. Math. Phys.} {\bf 61} (2002), 139--147.

\bibitem{gers} M. Gerstenhaber, On the deformation of rings and algebras, {\em Ann. of Math. (2)} {\bf 79} (1964), 59--103.

\bibitem{gls} L. Guo, H. Lang and Y. Sheng, Integration and geometrization of Rota-Baxter Lie algebras, {\em Adv. Math.} {\bf 387} (2021), 107834.

\bibitem{gqwz} S. Guo, Y. Qin, K. Wang and G. Zhou, Deformations and cohomology theory of Rota-Baxter 3-Lie algebras of arbitrary weights, {\em J. Geom. Phys.} {\bf 183} (2023), 104704.

\bibitem{hmns} F. Harrathi, S. Mabrouk, O. Ncib and S. Silvestrov, Malcev Yang-Baxter equation, weighted $\mathcal{O}$-operators on Malcev algebras and post-Malcev algebras, {\em Hacet. J. Math. Stat.} {\bf 52} (2023), no. 5, 1151--1171.

\bibitem{hsz} S. Hou, Y. Sheng and Y. Zhou, 3-post-Lie algebras and relative Rota-Baxter operators of nonzero weight on 3-Lie algebras, {\em J. Algebra} {\bf 615} (2023), 103--129.

\bibitem{kuper} B. A. Kupershmidt, What a classical $r$-matrix really is, {\em J. Nonlinear Math. Phys.} {\bf 6} (1999), 448--488.

\bibitem{lst} A. Lazarev, Y. Sheng and R. Tang, Deformations and homotopy theory of relative Rota-Baxter Lie algebras, {\em Comm. Math. Phys.} {\bf 383} (2021), 595--631.

\bibitem{loday-une} J.-L. Loday, Une version non commutative des alg\`{e}bres de Lie: les alg\`{e}bres de Leibniz, {\em Enseign. Math. (2)} {\bf 39} (1993), no. 3--4, 269--293.

\bibitem{loday-pira} J.-L. Loday and T. Pirashvili, Universal enveloping algebras of Leibniz algebras and (co)homology, {\em Math. Ann.} {\bf 296} (1993), no. 1, 139--158.

\bibitem{rota} G.-C. Rota, Baxter algebras and combinatorial identities. I, II, {\em Bull. Amer. Math. Soc.} {\bf 75} (1969), 325--329; ibid. 75 (1969), 330--334.

\bibitem{sun-guo} Q. Sun and S. Guo, Leibniz-dendriform bialgebras and relative Rota-Baxter operators, arXiv:2510.16826 (2025).

\bibitem{sun-jing} Q. Sun and N. Jing, $\mathcal{O}$-operators and related structures on Leibniz algebras, {\em Comm. Algebra} {\bf 51} (2023), no. 5, 2199--2216.

\bibitem{tang} R. Tang, C. Bai, L. Guo and Y. Sheng, Deformations and their controlling cohomologies of $\mathcal{O}$-operators, {\em Comm. Math. Phys.} {\bf 368} (2019), no. 2, 665--700.

\bibitem{sheng-pl} R. Tang and Y. Sheng, Leibniz bialgebras, relative Rota-Baxter operators, and the classical Leibniz Yang-Baxter equation, {\em J. Noncommut. Geom.} {\bf 16} (2022), 1179--1211.

\bibitem{sheng-defr} R. Tang, Y. Sheng and Y. Zhou, Deformations of relative Rota-Baxter operators on Leibniz algebras, {\em Int. J. Geom. Methods Mod. Phys.} {\bf 17} (2020), no. 12, 2050174.

\bibitem{tg-lya} W. Teng and S. Guo, Modified Rota-Baxter Lie-Yamaguti algebras, arXiv:2401.17726 (2024).

\bibitem{teng-hlts} W. Teng, Weighted $\mathcal{O}$-operators on Hom-Lie triple systems, {\em Complex Anal. Oper. Theory} {\bf 20} (2026), doi:10.1007/s11785-026-01908-6.

\bibitem{uchino} K. Uchino, Quantum analogy of Poisson geometry, related dendriform algebras and Rota-Baxter operators, {\em Lett. Math. Phys.} {\bf 85} (2008), no. 2--3, 91--109.

\bibitem{vallette} B. Vallette, Homology of generalized partition posets, {\em J. Pure Appl. Algebra} {\bf 208} (2007), 699--725.

\bibitem{voro} Th. Voronov, Higher derived brackets and homotopy algebras, {\em J. Pure Appl. Algebra} {\bf 202} (2005), no. 1--3, 133--153.

\end{thebibliography}
\end{document}